\documentclass[10pt]{amsart}

\usepackage{amsfonts, amsmath, amsthm, amssymb, epsfig, bm, graphics,
  psfrag, latexsym, chngcntr, color, yfonts, mathtools, phonetic,
  pigpen, subfig, longtable}
\usepackage[all,cmtip]{xy}
\usepackage{tikz} \usetikzlibrary{matrix,arrows}
\definecolor{darkblue}{rgb}{0,0,0.4}
\usepackage[colorlinks=true, citecolor=darkblue, filecolor=darkblue, linkcolor=darkblue, urlcolor=darkblue]{hyperref}
\usepackage[all]{hypcap}

\newtheorem{thm}{Theorem}[section]         
\newtheorem{lem}[thm]{Lemma}

\newtheorem{conj}[thm]{Conjecture}  
\counterwithin{figure}{section}

\newcommand{\R}{\mathbb{R}}
\newcommand{\Z}{\mathbb{Z}}
\newcommand{\F}{\mathbb{F}}
\newcommand{\C}{\mathbb{C}}
\newcommand{\T}{\mathbb{T}}

\newcommand{\mc}{\mathcal}
\newcommand{\mb}{\mathbb}
\newcommand{\mf}{\mathfrak}

\newcommand{\wh}{\widehat}
\newcommand{\wt}{\widetilde}
\newcommand{\ol}{\overline}

\newcommand{\del}{\partial}
\newcommand{\sbs}{\subset}
\newcommand{\sbseq}{\subseteq}
\newcommand{\sm}{\setminus}

\newcommand{\Si}{\Sigma}
\newcommand{\si}{\sigma}
\newcommand{\al}{\alpha}
\newcommand{\be}{\beta}

\newcommand{\vphi}{\varphi}

\newcommand{\Zoltan}{Zolt\'{a}n}
\newcommand{\Szabo}{Szab\'{o}}
\newcommand{\Ozsvath}{Ozsv\'{a}th}
\newcommand{\Andras}{Andr\'{a}s}
\newcommand{\Juhasz}{Juh\'{a}sz}

\newcommand{\ith}{^{\text{th}}}
\newcommand{\Aut}{\operatorname{Aut}}

\newcommand{\Mor}{\operatorname{Mor}}

\newcommand{\Ob}{\operatorname{Ob}}
\newcommand{\Id}{\operatorname{Id}}

\newcommand{\Sym}{\text{Sym}}
\newcommand{\nbd}{\text{nbd}}
\newcommand{\Diff}{\text{Diff}}
\newcommand{\MCG}{\text{MCG}}

\newcommand{\CF}{\mathit{CF}}

\newcommand{\CFL}{\mathit{CFL}}
\newcommand{\HFL}{\mathit{HFL}}
\newcommand{\GL}{\mathit{GL}}

\begin{document}
\title[Moving basepoints and the induced automorphisms of link Floer homology]{Moving basepoints and the induced automorphisms\\ of link Floer homology}
\author{Sucharit Sarkar}
\address{Department of Mathematics, Columbia University, New York, NY 10027}
\email{\href{mailto:sucharit.sarkar@gmail.com}{sucharit.sarkar@gmail.com}}

\subjclass[2010]{\href{http://www.ams.org/mathscinet/search/mscdoc.html?code=57M25,57M27,57R58}{57M25,
57M27, 57R58}}
\keywords{Link Floer homology; basepoint; mapping class group action; grid diagram}

\date{}

\begin{abstract}
  Given an $l$-component pointed oriented link $(L,p)$ in an oriented
  three-manifold $Y$, one can construct its link Floer chain complex
  $\CFL(Y,L,p)$ over the polynomial ring $\F_2[U_1,\ldots,U_l]$. Moving
  the basepoint $p_i\in L_i$ once around the link component $L_i$ induces
  an automorphism of $\CFL(Y,L,p)$. In this paper, we study an
  automorphism (a possibly different one) of $\CFL(Y,L,p)$ defined
  explicitly in terms of holomorphic disks; for links in $S^3$, we
  show that these two automorphisms are the same.
\end{abstract}

\maketitle

\section{Introduction}
Heegaard Floer theory is a collection of invariants, originally
defined for pointed oriented closed
three-manifolds \cite{POZSz,POZSzapplications}, and subsequently
extended for pointed oriented knots \cite{POZSzknotinvariants,JR} and
pointed oriented links \cite{POZSzlinkinvariants} in oriented
three-manifolds. For each of these objects, the theory comes in
several variants; in each variant, one constructs a chain complex in
the graded homotopy category over some graded ring; furthermore, for
each object and in each variant, the mapping class group of the object
acts on the chain complex.

In this paper, we will work with links. Initially, we will study links
in arbitrary oriented three-manifolds; later on, we will concentrate
on links in $S^3$. For the experts only\footnote{Non-experts should
skip to the next paragraph.}, let us specify the version of link Floer
theory that we will study: it is the associated graded object of the
minus version of the fully filtered theory over the base ring
$\F_2[U_1,\ldots,U_l]$. We choose to work with the associated graded
object and not the (more general) fully filtered theory because the
former object is easier to study and also because the proof
of \hyperref[thm:main]{Theorem \ref*{thm:main}} is simpler. For
simplicity again, we only work with the minus version $\CFL^{-}$ and
not all three of the standard versions $\CFL^{-},\CFL^{+},\CFL^{\infty}$
(the general object that we could have studied being the inclusion
$\iota\colon\CFL^{-}\hookrightarrow \CFL^{\infty}$); however, we can
construct $\CFL^{\infty}$ and $\CFL^{+}$ formally from $\CFL^{-}$ as
$\CFL^{-}\otimes_{\F_2[U_1,\ldots,U_l]}\F_2[U_1,U_1^{-1},\ldots,U_l,U_l^{-1}]$
and the mapping cone of
$\CFL^{-}\rightarrow \CFL^{-}\otimes_{\F_2[U_1,\ldots,U_l]}\F_2[U_1,U_1^{-1},\ldots,U_l,U_l^{-1}]$,
respectively.
We work over the ring $\F_2[U_1,\ldots,U_l]$ and not the
more universal ring $\Z[U_1,\ldots,U_l]$ because most of the variants
of link Floer theory are only defined over the base field $\F_2$ (and
again, because working over $\F_2$ is simpler). Therefore, henceforth
whenever we say $\CF$ or $\CFL$, we mean the associated graded object of the
minus version of the link Floer chain complex over
$\F_2[U_1,\ldots,U_l]$.

As is standard in Heegaard Floer theory, a pointed link $(L,p)$ in a
thee-manifold $Y$ is described by a Heegaard diagram $\mc{H}$, and a
chain complex $\CF_{(\mc{H},J_s)}$ is defined which depends on
$\mc{H}$ and some additional data $J_s$; the link-invariant
$\CFL(Y,L,p)$ can be obtained from $\CF_{(\mc{H},J_s)}$ using
naturality. We will define certain link-invariant maps $\Phi_i$ and
$\Psi_i$ from $\CF_{(\mc{H},J_s)}$ to itself using counts of certain
holomorphic disks. Our main theorem is the following.

\begin{thm}\label{thm:main}
  Let $\mc{H}=(\Si,\al,\be,z,w)$ be a Heegaard diagram representing an
  $l$-component pointed link $(L,p)$ in $S^3$; then for all $1\leq
  i\leq l$, the automorphism $\Id+\Psi_{i}\Phi_{i}$ in
  $\Aut_{K(\mc{A}_l)}(\CF_{(\mc{H},J_s)})$ induces the automorphism
  $\rho(\si_i)$ in $\Aut_{N(K(\mc{A}_l))}(\CFL(S^3,L,p))$, where
  $\si_i\in\MCG(S^3,L,p)$ is the positive Dehn twist along $i\ith$
  link component $L_i$ and $\rho(\si_i)$ is its induced automorphism
  on $\CFL(S^3,L,p)$.
\end{thm}

In \hyperref[sec:basics]{Section \ref*{sec:basics}}, we will give a
quick tour of the relevant areas of Heegaard Floer theory; in
\hyperref[sec:mcg]{Section \ref*{sec:mcg}}, we will talk about the
mapping class group action, and a specific mapping class group element
$\si_i$, the positive Dehn twist around the $i\ith$ link component
$L_i$, which corresponds to moving the basepoint $p_i\in L_i$ once
around; in \hyperref[sec:candidate]{Section \ref*{sec:candidate}}, we
will define the maps $\Phi_i$ and $\Psi_i$ and prove that the map
$\Id+\Phi_i\Psi_i$ induces a well-defined automorphism of
$\CFL(Y,L,p)$; in \hyperref[sec:grids]{Section \ref*{sec:grids}}, we
will use grid diagrams to prove \hyperref[thm:main]{Theorem
  \ref*{thm:main}}; and finally in \hyperref[sec:computations]{Section
  \ref*{sec:computations}}, we will compute the automorphism
$\Id+\Phi_1\Psi_1$ for all the $85$ prime knots up to nine crossings,
and see that it is non-trivial (as in, not the identity) more often
than not.

\subsection*{Acknowledgment}
The author was supported by the Clay Postdoctoral Fellowship when this
paper was written. He would like to thank Elisenda
Grigsby, \Andras{} \Juhasz{}, Robert Lipshitz, Peter \Ozsvath{}
and \Zoltan{} \Szabo{} for many helpful discussions.

\section{Heegaard Floer basics}\label{sec:basics}
A \emph{pointed link} is a link with a basepoint in each
component. Let $L\sbs Y$ be an oriented $l$-component pointed link
inside a closed oriented three-manifold $Y$; let $L_i$ be the $i\ith$
component, and let $p_i\in L_i$ be the basepoint in the $i\ith$
component. We say that $\mc{H}=(\Si,\al,\be,z,w)$ is a \emph{Heegaard
diagram} for $L$ if there exists a self-indexing Morse function
$f\colon Y\rightarrow\R$, equipped with a gradient-like flow, such that:
$\Si=f^{-1}(\frac{3}{2})$ is a surface of genus $g$;
$z=(z_1,\ldots,z_l)$, $w=(w_1,\ldots,w_l)$, and there is an $l$-tuple
$k=(k_1,\ldots,k_l)$, such that $z_i$ is a collection of $k_i$
markings $z_{i,1},\ldots,z_{i,k_i}$ in $\Si$ and $w_i$ is also a
collection of $k_i$ markings $w_{i,1},\ldots,w_{i,k_i}$ in $\Si$; for
each $i\in\{1,\ldots,l\}$, one of the $w_i$-markings, say $w_{i,s_i}$,
is designated \emph{special}; $f$ has $|k|$ index-zero critical points
and $|k|$ index-three critical points; $\al$ is the intersection of
$\Si$ and the stable manifold of the index-one critical points; $\be$
is the intersection of $\Si$ and the unstable manifold of the
index-two critical points; $L_i$ is the union of the flowlines through
the $z_i$-markings and the reversed flowlines through the
$w_i$-markings; the basepoint $p_i\in L_i$ is the special
$w_i$-marking $w_{i,s_i}$; the Heegaard diagram is also assumed to be
admissible \cite[Definition 3.5]{POZSzlinkinvariants}.

Let $n=g+|k|-1$, and let $\mf{j}$ be a complex structure on $\Si$; $J_s$
is a path of nearly symmetric almost complex structures on the
symmetric product $\Sym^n(\Si)$, which is a generic perturbation of
the constant path $\Sym^n(\mf{j})$ \cite[Definition 3.1]{POZSz};
$\T_{\al}=\{x\in\Sym^n(\Si)\mid\text{all the coordinates of }x\text{
  lie on }\al\}$ is a totally real half-dimensional torus; $\T_{\be}$
is defined similarly; the marking $z_{i,j}$ gives rise to the divisor
$Z_{i,j}=\{x\in\Sym^n(\Si)\mid\text{at least one of the coordinates of
}x\text{ is }z_{i,j}\}$; $W_{i,j}$ is defined similarly; let
$Z_i=\sum_jZ_{i,j}$, $W_i=\sum_jW_{i,j}$, $Z=\sum_i Z_i$ and
$W=\sum_iW_i$.

Given $x,y\in\T_{\al}\cap\T_{\be}$, $\pi_2(x,y)$ is the set of all
\emph{Whitney disks} joining $x$ to $y$, or in other words, the set of
all homotopy classes of maps $(\{z\in\C\mid|z|\leq
1\},\{ie^{i\theta}\mid
-\frac{\pi}{2}\leq\theta\leq\frac{\pi}{2}\},\{e^{i\theta}\mid
-\frac{\pi}{2}\leq\theta\leq\frac{\pi}{2}\},i,-i)
\rightarrow(\Sym^n(\Si),\T_{\al},\T_{\be},x,y)$; given a Whitney disk
$\vphi\in\pi_2(x,y)$, $\wh{\mc{M}}_{J_s}(\vphi)$
is the \emph{unparametrized moduli space} of such maps. Elements of
$\T_{\al}\cap\T_{\be}$ carry a \emph{Maslov grading} $M$ \cite[Theorem
  7.1]{POZSz4manifolds} and $l$ \emph{Alexander gradings} $A_i$
\cite[Subsection 8.1]{POZSzlinkinvariants}, such that whenever
$\vphi\in\pi_2(x,y)$, $M(y)-M(x)=\mu(\vphi)-2\vphi\cdot W$ and
$A_i(x)-A_i(y)=\vphi\cdot (Z_i-W_i)$, where the \emph{Maslov index}
$\mu(\vphi)$ is the expected dimension the moduli space
$\mc{M}_{J_s}(\vphi)$. Let $\mb{P}$ be the $(l+1)$-graded polynomial
ring generated over $\F_2$ by the variables $U_{i,j}$ for
$i\in\{1,\ldots,l\},j\in\{1,\ldots,k_i\}$, where the
$(M,A_1,\ldots,A_l)$ grading of $U_{i,j}$ is
$(-2,-\delta_{1i},\ldots,-\delta_{li})$\footnote{Throughout the paper,
$\delta$ denotes the Kronecker delta function.}. The chain complex
$\CF_{(\mc{H},J_s)}$ is the $(l+1)$-graded
$\F_2[U_1,\ldots,U_l]$-module freely generated over $\mb{P}$ by
$\T_{\al}\cap\T_{\be}$, where the $U_i$-action is multiplication by
$U_{i,s_i}$; the boundary map is an $U_{i,j}$-equivariant
$(-1,0,\ldots,0)$-graded map, and for $x\in\T_{\al}\cap\T_{\be}$, it
is given by
$$\del x=\sum_{y\in\T_{\al}\cap
  T_{\be}}y\sum_{\substack{\vphi\in\pi_2(x,y)\\\vphi\cdot
  Z=0\\\mu(\vphi)=1}}|\wh{\mc{M}}_{J_s}(\vphi)|\prod_{\imath,\jmath}U^{\vphi\cdot
  W_{\imath,\jmath}}_{\imath,\jmath}.$$

Given a small Abelian category $\mc{C}$, let $K(\mc{C})$ be the
homotopy category of chain complexes whose objects are chain complexes
in $\mc{C}$, and whose morphisms are chain maps up to chain homotopy.
Given a small category $\mc{C}$ and a group $G$, let $\mc{C}_G$ be the
category whose objects are two-tuples $(A,f)$, where
$A\in\Ob_{\mc{C}}$ and $f\colon G\rightarrow\Aut_{\mc{C}}(A)$, and the set
of morphisms $\Mor_{\mc{C}_G}((A_1,f_1),(A_2,f_2))$ is the subset of
$\Mor_{\mc{C}}(A_1,A_2)$ consisting of the ones that are
$G$-equivariant.  Given a small category $\mc{C}$, let $N(\mc{C})$ be
the category whose objects are three-tuples $(I,\mf{o}_I,\mf{f}_I)$,
where $I$ is a set, $\mf{o}_I$ is a map from $I$ to $\Ob_{\mc{C}}$ and
$\mf{f}_I$ is a map from $I\times I$ to $\Mor_{\mc{C}}$, such that:
$\mf{f}_I(i,i')\in\Mor_{\mc{C}}(\mf{o}_I(i),\mf{o}_I(i'))$ for all
$i,i'\in I$; $\mf{f}_I(i,i)=\Id_{\mf{o}_I(i)}$ for all $i\in I$; and
$\mf{f}_I(i',i'')\mf{f}_I(i,i')=\mf{f}_I(i,i'')$ for all $i,i',i''\in
I$. A morphism from $(I,\mf{o}_I,\mf{f}_I)$ to $(J,\mf{o}_J,\mf{f}_J)$
is a map $\phi$ from $I\times J$ to $\Mor_{\mc{C}}$, such that:
$\phi(i,j)\in \Mor_{\mc{C}}(\mf{o}_I(i),\mf{o}_J(j))$ for all $i\in I$
and all $j\in J$; $\phi(i',j)\mf{f}_I(i,i')=\phi(i,j)$ for all
$i,i'\in I$ and all $j\in J$; and $\mf{f}_J(j,j')\phi(i,j)=\phi(i,j')$
for all $i\in I$ and all $j,j'\in J$.

If $\mc{A}_l$ is the category of $(A_1,\ldots,A_l)$-graded
$\F_2[U_1,\ldots,U_l]$-modules, then $\CF_{(\mc{H},J_s)}$ is an object
of $K(\mc{A}_l)$.  If $\mc{H}'=(\Si',\al',\be',z',w')$ is another
Heegaard diagram for the same pointed link $(L,p)$ and if $J'_s$ is a
path of nearly symmetric almost complex structures on
$\Sym^{n'}(\Si')$, then by naturality \cite[Theorem
2.1]{POZSz4manifolds} \cite[Section 6]{POAScontactsurgeries}
\cite[Subsection 5.2]{AJcobordism} \cite{AJPODT}, there is an
$U_i$-equivariant chain map from $\CF_{(\mc{H},J_s)}$ to
$\CF_{(\mc{H}',J'_s)}$.  It can be checked that this map is
well-defined up to $U_i$-equivariant chain homotopy, or in other
words, this map induces a well-defined morphism
$F_{(\mc{H},J_s),(\mc{H}',J'_s)}$ in $K(\mc{A}_l)$.  Therefore, given
a pointed link $(L,p)$ in $Y$, we get a well-defined object
$\CFL(Y,L,p)$ in $N(K(\mc{A}_l))$, where the indexing set $I$ is the
set of all ordered pairs $(\mc{H},J_s)$, where $\mc{H}$ is a Heegaard
diagram for $L$ and $J_s$ is a path of nearly symmetric almost complex
structures on the symmetric product, $\mf{o}_I((\mc{H},J_s))$ is the
chain complex $\CF_{(\mc{H},J_s)}$, and
$\mf{f}_I((\mc{H},J_s),(\mc{H}',J'_s))$ is the morphism
$F_{(\mc{H},J_s),(\mc{H}',J'_s)}$.

Let $\mc{B}_l$ be the category of $(M,A_1,\ldots,A_l)$-graded
$\F_2[U_1,\ldots, U_l]$-modules, and let $\mc{C}_l$ be the category of
$(M,A_1,\ldots,A_l)$-graded $\F_2$-modules. By taking homology, we get
a pointed link invariant object $\HFL(Y,L,p)\allowbreak
=H_*(\CFL(Y,L,p))$ in $N(\mc{B}_l)$; after putting $U_i=0$ for all
$1\leq i\leq l$, and then taking homology, we get a pointed link
invariant object $\wh{\HFL}(Y,L,p)=H_*(\CFL(Y,L,p)/\{U_i=0\})$ in
$N(\mc{C}_l)$. We sometimes need the shift functors in $K(\mc{A}_l)$,
$\mc{B}_l$ and $\mc{C}_l$ (and the induced shift functors in
$N(K(\mc{A}_l))$, $N(\mc{B}_l)$ and $N(\mc{C}_l)$). Let
$[m,a_1,\ldots,a_l]$ be the shift functor in any of these categories
that decreases the $(M,A_1,\ldots,A_l)$-grading by
$(m,a_1,\ldots,a_l)$.

\section{Mapping class group actions}\label{sec:mcg}
Diffeomorphisms will always be orientation preserving\footnote{In
  particular, unless otherwise mentioned, a diffeomorphism of the pair
  $(Y,L)$ is a diffeomorphism that preserves the orientations of both
  $Y$ and $L$.}, and the
\emph{mapping class group} $\MCG$ is the $\pi_0$ of the space of all
(orientation preserving) self-diffeomorphisms. There exists an
well-defined map
$\rho\colon\MCG(Y,L,p)\rightarrow\Aut_{N(K(\mc{A}_l))}(\CFL(Y,L,p))$,
defined in \cite{POZSz4manifolds} \cite[Definition
6.5]{POAScontactsurgeries} \cite[Corollary 5.20]{AJcobordism}
\cite{AJPODT} as follows. Let $\si\in\MCG(Y,L,p)$; assume that $\si$
comes from $\wt{\si}\in\Diff(Y,L,p)$; let $\wt{\si}$ also denote the
induced automorphisms of the set of all Heegaard diagrams for
$(Y,L,p)$, their symmetric products, and the space of all paths of
nearly symmetric almost complex structures on the symmetric
products. Therefore, we have an $U_{i,j}$-equivariant chain map from
$\CF_{(\wt{\si}^{-1}(\mc{H}),\wt{\si}^{-1}(J_s))}$ to
$\CF_{(\mc{H},J_s)}$ which sends
$x\in\wt{\si}^{-1}(\T_{\al})\cap\wt{\si}^{-1}(\T_{\be})$ to
$\wt{\si}(x)\in\T_{\al}\cap\T_{\be}$. This induces a well-defined
morphism $f_{\si}$ from
$\CF_{(\wt{\si}^{-1}(\mc{H}),\wt{\si}^{-1}(J_s))}$ to
$\CF_{(\mc{H},J_s)}$ in $K(\mc{A}_l)$, and hence a well-defined
automorphism
$f_{\si}F_{(\mc{H},J_s),(\wt{\si}^{-1}(\mc{H}),\wt{\si}^{-1}(J_s))}$
in $\Aut_{K(\mc{A}_l)}(\CF_{(\mc{H},J_s)})$. It turns out that the
maps of the form $F_{(\mc{H},J_s),(\mc{H}',J'_s)}$ commute with such
automorphisms. Therefore, we can treat $\CF_{(\mc{H},J_s)}$ as an
object in $K(\mc{A}_l)_{\MCG(Y,L,p)}$ and the maps
$F_{(\mc{H}_,J_s),(\mc{H},J_s)}$ as morphisms in this
category. Therefore, $\CFL(Y,L,p)$ can be thought of as a pointed link
invariant object in $N(K(\mc{A}_l))_{\MCG(Y,L,p)}$.

Let $T(L)=\prod_i L_i$ be the pointed $l$-dimensional torus with the
basepoint $p=(p_1,\ldots,p_l)$. Since $L$ is oriented, $\pi_1(T(L),p)$
is canonically isomorphic to $\Z^l$. We have a fiber bundle
$$\begin{tikzpicture}
\matrix(m)[matrix of math nodes,
row sep=2.5em, column sep=2em,
text height=2ex, text depth=0.25ex]
{(\Diff(Y,L,p),\Id)&(\Diff(Y,L),\Id)\\
&(T(L),p)\\};
\path[right hook->](m-1-1) edge (m-1-2);
\path[->](m-1-2) edge (m-2-2);
\end{tikzpicture}$$
which gives rise to a long exact sequence
$$\Z^l=\pi_1(T(L),p)\rightarrow \MCG(Y,L,p)= \pi_0(\Diff(Y,L,p))\rightarrow \MCG(Y,L)
=\pi_0(\Diff(Y,L))\rightarrow \{0\}.$$ 
Let $\si_i\in\MCG(Y,L,p)$ be the
image of the $i\ith$ unit vector in $\Z^l$; we call $\si_i$
the \emph{positive Dehn twist around $L_i$}. Then there is an action of
$\MCG(Y,L)$ on $\CFL(Y,L,p)$, which is well-defined up to the
$l$ automorphisms $\rho(\si_i)$. In this paper, we will try to
understand these $l$ automorphisms.

\begin{figure}
\psfrag{a}{$\al$}
\psfrag{a'}{$\al'$}
\psfrag{b}{$\be$}
\psfrag{b'}{$\be'$}
\psfrag{z}{$z_{i,1}$}
\psfrag{w}{$w_{i,1}$}
\begin{center}
\includegraphics[width=0.5\textwidth]{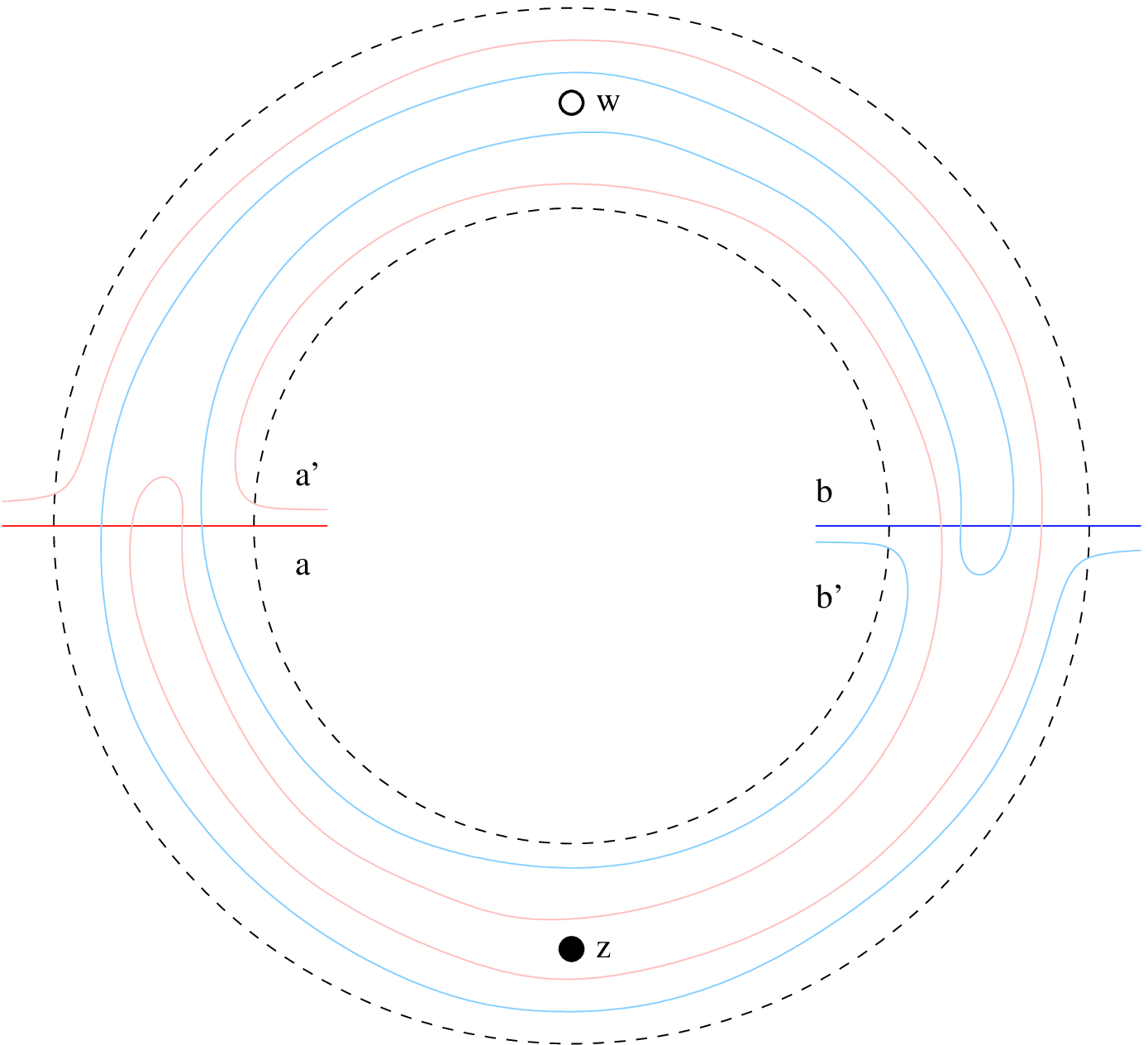}
\end{center}
\caption{The annular neighborhood $\nbd(C)$.}\label{fig:annular}
\end{figure}

Let us first describe a way to view $\rho(\si_i)$ as a composition of
two triangle maps. Let $\mc{H}_{\al\be}=(\Si,\al,\be,z,w)$ be a
Heegaard diagram for $L$ such that the $i\ith$ link component $L_i$
contains exactly one $w$ marking and exactly one $z$ marking. By
stabilizing twice if necessary, we can assume that there is an
oriented arc joining $z_{i,1}$ to $w_{i,1}$ which is disjoint from
$\al$ and intersects $\be$ transversely at a point, and there is an
oriented arc joining $w_{i,1}$ to $z_{i,1}$ which is disjoint from
$\be$ and intersects $\al$ transversely at a point, and the union of
these two arcs is an oriented embedded circle $C$ on the Heegaard
surface $\Si$. A regular neighborhood $\nbd(C)$ of $C$ is shown in
\hyperref[fig:annular]{Figure \ref*{fig:annular}}.

Let $\del(\nbd(C))=C_1-C_2$, where each $C_i$ is oriented parallel to
$C$. Let $\si\in\Diff(\Si,z,w)$ be the composition of a positive Dehn
twist along $C_1$ and a negative Dehn twist along $C_2$. Let $J_s$ be
a path of nearly symmetric almost complex structures on $\Sym^n(\Si)$;
and let $J_{s,t}$ be a path in the space of paths of nearly symmetric
almost complex structures joining $\si^{-1}(J_s)$ to $J_s$. Let $\al'$
be obtained by first perturbing $\al$ and then applying $\si^{-1}$; $\be'$
is defined similarly. The multicurves $\al$, $\be$, $\al'$ and $\be'$
in $\nbd(C)$ are represented
in \hyperref[fig:annular]{Figure \ref*{fig:annular}} by the red, blue,
pink and light blue curves, respectively. Let
$\mc{H}_{\al\be'}=(\Si,\al,\be',z,w)$ and
$\mc{H}_{\al'\be'}=(\Si,\al',\be',z,w)$.

The $\MCG(Y,L,p)$ element induced by $\si$ is $\si_i$. Therefore, the
automorphism $\rho(\si_i)$ acts by mapping
$x\in\T_{\al'}\cap\T_{\be'}$ in
$\CF_{(\mc{H}_{\al'\be'},\si^{-1}(J_s))}$ to
$\si(x)\in\T_{\al}\cap\T_{\be}$ in
$\CF_{(\mc{H}_{\al\be},J_s)}$. The naturality map
$F_{(\mc{H}_{\al\be},J_s),(\mc{H}_{\al'\be'},\si^{-1}(J_s))}$, by
naturality, is the composition
$F_{(\mc{H}_{\al'\be'},J_s),(\mc{H}_{\al'\be'},\si^{-1}(J_s))}
F_{(\mc{H}_{\al\be'},J_s),(\mc{H}_{\al'\be'},J_s)}
F_{(\mc{H}_{\al\be},J_s),(\mc{H}_{\al\be'},J_s)}$. The first two
maps are the triangle maps \cite[Equation 21]{POZSz}, and it is easily
verified that the relevant triple Heegaard diagrams are admissible.
The third map is induced by the path $J_{s,t}$ \cite[Equation
14]{POZSz}; however, if we assume that $\wh{M}_{J_{s,t}}(\vphi)$ is
empty for all $t\in[0,1]$ and for all Whitney disks $\vphi$ with
$\mu(\vphi)\leq 0$, for example by assuming that $\mc{H}_{\al\be}$ is
a nice Heegaard diagram \cite[Definition 3.1]{SSJW}, then the third
map is the identity map. Therefore, the automorphism $\rho(\si_i)$ can
be thought of as a composition of two triangle maps.

\section{A candidate}\label{sec:candidate}

Following the notations
from \hyperref[sec:basics]{Section \ref*{sec:basics}}, let
$\mc{H}=(\Si,\al,\be,z,w)$ be a Heegaard diagram for a pointed link
$(L,p)$ in $Y$, and let $J_s$ be a path of nearly symmetric almost
complex structures on $\Sym^n(\Si)$. For each $1\leq i\leq l$, let
$a_i,b_i\colon\{1,\ldots,k_i\}\rightarrow\{1,\ldots,k_i\}$ be the two
bijections such that for every $j$, $z_{i,j}$ and $w_{i,a_i(j)}$ lie
in the same component of $\Si\sm\al$ and $z_{i,j}$ and $w_{i,b_i(j)}$
lie in the same component of $\Si\sm\be$. For each $(i,j)$, let us
define two $\mb{P}$-module maps $\Psi_{i,j}$ and $\Phi_{i,j}$
from $\CF_{(\mc{H},J_s)}$ to
$\CF_{(\mc{H},J_s)}[-1,-\delta_{1i},\ldots,-\delta_{li}]$ and
$\CF_{(\mc{H},J_s)}[1,\delta_{1i},\ldots,\delta_{li}]$
respectively, as follows: for $x\in\T_{\al}\cap\T_{\be}$,
\begin{align*}
\Psi_{{i,j}}(x)=&\sum_{y\in\T_{\al}\cap
  T_{\be}}y\sum_{\substack{\vphi\in\pi_2(x,y)\\\vphi\cdot
  Z=\vphi\cdot
  Z_{i,j}=1\\\mu(\vphi)=1}}|\wh{\mc{M}}_{J_s}(\vphi)|
\prod_{\imath,\jmath}U^{\vphi\cdot
  W_{\imath,\jmath}}_{\imath,\jmath}\text{ and}\\
\Phi_{{i,j}}(x)=&\sum_{y\in\T_{\al}\cap
  T_{\be}}y\sum_{\substack{\vphi\in\pi_2(x,y)\\\vphi\cdot
  Z=0\\\mu(\vphi)=1}}(\vphi\cdot
  W_{i,j})|\wh{\mc{M}}_{J_s}(\vphi)|U^{-1}_{i,j} \prod_{\imath,\jmath}U^{\vphi\cdot
  W_{\imath,\jmath}}_{\imath,\jmath}.
\end{align*}

\begin{lem}
For every $x\in\T_{\al}\cap\T_{\be}$, and for every $(i,j)$, the
commutator
$[\del\colon\Psi_{i,j}]=U_{i,a_i(j)}+U_{i,b_i(j)}$
and the commutator $[\del\colon\Phi_{{i,j}}]=0$.
\end{lem}

\begin{proof}
  Let us first set up a few extra notations. For any $x\in\T_{\al}$,   let $\pi_2^{\al}(x)$ be the set of all Whitney disks with boundary   lying in $\T_{\al}$, or in other words, the set of all homotopy   classes of maps $(\{z\in\C\mid|z|\leq   1\},\{z\in\C\mid|z|=1\},i)\rightarrow (\Sym^n(\Si),\T_{\al},x)$;   given $\vphi^{\al}\in\pi_2^{\al}(x)$,   $\wh{\mc{N}}^{\al}_{J_s}(\vphi)$ is the unparametrized moduli space   of such maps. The Whitney disks $\pi_2^{\be}(x)$ and the moduli   spaces $\wh{\mc{N}}^{\be}_{J_s}$ are defined similarly. For every   $x\in\T_{\al}$ and for every $(i,j)$ with $1\leq i\leq l$, $1\leq   j\leq k_i$, there exists a unique Whitney disk   $\vphi^{\al}_{x,i,j}\in\pi_2^{\al}(x)$ such that   $\mu(\vphi^{\al}_{x,i,j})=2$ and $\vphi^{\al}_{x,i,j}\cdot   Z_{\imath,\jmath}$ is $1$ if $(\imath,\jmath)=(i,j)$, and is $0$ otherwise; also,   $\vphi^{\al}_{x,i,j}\cdot W_{\imath,\jmath}$ is $1$ if $(\imath,\jmath)=(i,a_i(j))$,   and is $0$ otherwise; furthermore, if $\vphi^{\al}\in\pi_2^{\al}(x)$   with $\mu(\vphi^{\al})=2$, then   $\wh{\mc{N}}^{\al}_{J_s}(\vphi^{\al})$ has an odd number of points   if and only if $|k|>1$ and $\vphi^{\al}=\vphi^{\al}_{x,i,j}$ for   some $(i,j)$ \cite[Theorem 5.5]{POZSzlinkinvariants}. Similarly, for   every $x\in\T_{\be}$ and for every $(i,j)$, there exists a unique   Whitney disk $\vphi^{\be}_{x,i,j}\in\pi_2^{\be}(x)$ such that   $\mu(\vphi^{\be}_{x,i,j})=2$ and $\vphi^{\be}_{x,i,j}\cdot   Z_{\imath,\jmath}$ is $1$ if $(\imath,\jmath)=(i,j)$, and is $0$ otherwise; also,   $\vphi^{\be}_{x,i,j}\cdot W_{\imath,\jmath}$ is $1$ if $(\imath,\jmath)=(i,b_i(j))$,   and is $0$ otherwise; furthermore, if $\vphi^{\be}\in\pi_2^{\be}(x)$   with $\mu(\vphi^{\be})=2$, then   $\wh{\mc{N}}^{\be}_{J_s}(\vphi^{\be})$ has an odd number of points   if and only if $|k|>1$ and $\vphi^{\be}=\vphi^{\be}_{x,i,j}$ for   some $(i,j)$.

From Gromov compactification adapted to our present
settings \cite[Section 6]{POZSzlinkinvariants}, we know that for any
$\vphi\in\pi_2(x,y)$ with $\mu(\vphi)=2$, the number of broken
flowlines,
$\sum_{\vphi=\vphi_1*\vphi_2}|\wh{\mc{M}}_{J_s}(\vphi_1)\times\wh{\mc{M}}_{J_s}(\vphi_2)|$,
is even, unless $|k|>1$ and $\vphi=\vphi^{\al}_{x,i,j}$ or
$\vphi=\vphi^{\be}_{x,i,j}$ for some $(i,j)$, in which case, it is
odd. Therefore, for every $x\in\T_{\al}\cap\T_{\be}$, and for every $(i,j)$,
\begin{align*}
(\del\Psi_{{i,j}}+\Psi_{{i,j}}\del)(x)=&\sum_{y\in\T_{\al}\cap
  T_{\be}}y\sum_{\substack{\vphi\in\pi_2(x,y)\\\vphi\cdot
  Z=\vphi\cdot
  Z_{i,j}=1\\\mu(\vphi)=2}} 
\sum_{\substack{\vphi_1,\vphi_2\\\mu(\vphi_1)=\mu(\vphi_2)=1\\ \vphi=\vphi_1*\vphi_2}} 
  |\wh{\mc{M}}_{J_s}(\vphi_1)\times\wh{\mc{M}}_{J_s}(\vphi_2)|
\prod_{\imath,\jmath}U^{\vphi\cdot
  W_{\imath,\jmath}}_{\imath,\jmath}\\
=&(U_{i,a_i(j)}+U_{i,b_i(j)})x,\text{ and}\\
(\del\Phi_{{i,j}}+\Phi_{{i,j}}\del)(x)=&\sum_{y\in\T_{\al}\cap
  T_{\be}}y\sum_{\substack{\vphi\in\pi_2(x,y)\\\vphi\cdot
  Z=0\\\mu(\vphi)=2}}
\sum_{\substack{\vphi_1,\vphi_2\\\mu(\vphi_1)=\mu(\vphi_2)=1\\ \vphi=\vphi_1*\vphi_2}}
(\vphi_1+\vphi_2)\cdot
  W_{i,j}|\wh{\mc{M}}_{J_s}(\vphi_1)\times\wh{\mc{M}}_{J_s}(\vphi_2)|
  U^{-1}_{i,j} \prod_{\imath,\jmath}U^{\vphi\cdot W_{\imath,\jmath}}_{\imath,\jmath}\\
=&\sum_{y\in\T_{\al}\cap
  T_{\be}}y\sum_{\substack{\vphi\in\pi_2(x,y)\\\vphi\cdot
  Z=0\\\mu(\vphi)=2}}(\vphi\cdot W_{i,j})U_{i,j}^{-1}
\sum_{\substack{\vphi_1,\vphi_2\\\mu(\vphi_1)=\mu(\vphi_2)=1\\ \vphi=\vphi_1*\vphi_2}}
|\wh{\mc{M}}_{J_s}(\vphi_1)\times\wh{\mc{M}}_{J_s}(\vphi_2)|
  \prod_{\imath,\jmath}U^{\vphi\cdot W_{\imath,\jmath}}_{\imath,\jmath}\\
=&0\qedhere.
\end{align*}
\end{proof}

For $1\leq i\leq l$, let $\Psi_{i}=\sum_j\Psi_{{i,j}}$ and
let $\Phi_{i}=\Phi_{{i,s_i}}$

\begin{thm}\label{thm:invariance}
  For every $i$, the two $U_{i,j}$-equivariant chain maps $\Psi_{i}$
  and $\Phi_{i}$ induce link-invariant maps from $\CFL(Y,L,p)$ to
  $\CFL(Y,L,p)[-1,-\delta_{1i},\ldots,-\delta_{li}]$ and
  $\CFL(Y,L,p)[1,\delta_{1i},\ldots,\delta_{li}]$ respectively, in
  $N(K(\mc{A}_l))$.
\end{thm}

The proof is a consequence of the
following lemma.

\begin{lem}
  Let $\mc{H}'=(\Si',\al',\be',z',w')$ be another Heegaard diagram for
  $(L,p)$, and let $J'_s$ be a path of nearly symmetric almost complex
  structures on $\Sym^{n'}(\Si')$. If $\Psi'_{i}$ and
  $\Phi'_{i}$ denote the two chain maps on
  $\CF_{(\mc{H}',J'_s)}$, then for all $i$, both the maps
  $F_{(\mc{H},J_s),(\mc{H}',J'_s)}\Psi_{i}+\Psi'_{i}F_{(\mc{H},J_s),(\mc{H}',J'_s)}$
  and
  $F_{(\mc{H},J_s),(\mc{H}',J'_s)}\Phi_{i}+\Phi'_{i}F_{(\mc{H},J_s),(\mc{H}',J'_s)}$
  are chain homotopic to zero, where the chain homotopies are also
  $U_i$-equivariant.
\end{lem}

\begin{proof}
Let us first clarify a notational convention that we will follow for
the rest of the proof, and indeed, occasionally during the rest of the
paper. If an object in the Heegaard diagram $\mc{H}$
is denoted by some symbol $\textgoth{S}$, then the corresponding
object in the Heegaard diagram $\mc{H}'$ is denoted by
$\textgoth{S}'$. For example, the boundary map in $\mc{H}'$ is denoted
by $\del'$, the markings on the $i\ith$ link component $L_i$ are
$z'_{i,1},\ldots,z'_{i,k'_i},w'_{i,1},\ldots,w'_{i,k'_i}$, and the
basepoint $p_i$ is $w'_{i,s'_i}$.

Using \cite[Proposition 7.1]{POZSz} and \cite[Lemma 2.4]{CMPOSS}, we
can assume that we are in one of the following four cases.

\emph{\textbf{Case 1}: $\mc{H}'=\mc{H}$.} Following the proof of \cite[Theorem
6.1]{POZSz}, we assume that $\mf{j}'=\mf{j}$, and we choose a path
$J_{s,t}$ in the space of paths of nearly symmetric almost complex
structures joining $J_s$ to $J'_s$. Given a Whitney disk
$\vphi\in\pi_2(x,y)$, let $\mc{M}_{J_{s,t}}(\vphi)$ denote the moduli space of
holomorphic disks with time-dependent complex structure on the
target \cite[Equation 14]{POZSz}. The map
$F_{(\mc{H},J_s),(\mc{H},J'_s)}$ is defined as 
$$F_{(\mc{H},J_s),(\mc{H},J'_s)}(x)=\sum_{y\in\T_{\al}\cap
  T_{\be}}y\sum_{\substack{\vphi\in\pi_2(x,y)\\\vphi\cdot
  Z=0\\\mu(\vphi)=0}}|\mc{M}_{J_{s,t}}(\vphi)|
\prod_{\imath,\jmath}U^{\vphi\cdot
  W_{\imath,\jmath}}_{\imath,\jmath}.$$  
Define $H_{z_i},H_{w_i}\colon\CF_{(\mc{H},J_s)}\rightarrow
\CF_{(\mc{H},J'_s)}$\footnote{Technically there are some shift operators involved:
$H_{z_i}$ maps $\CF_{(\mc{H},J_s)}$ to $\CF_{(\mc{H},J'_s)}[0,-\delta_{1i},\ldots,-\delta_{li}]$ and
$H_{w_i}$ maps $\CF_{(\mc{H},J_s)}$ to
$\CF_{(\mc{H},J'_s)}[2,\delta_{1i},\ldots,\delta_{li}]$. However,
we will often suppress the degree shift information.} as
follows:
\begin{align*}
H_{z_i}(x)=&\sum_{y\in\T_{\al}\cap
  T_{\be}}y\sum_{\substack{\vphi\in\pi_2(x,y)\\\vphi\cdot Z=\vphi\cdot
  Z_i=1\\\mu(\vphi)=0}}|\mc{M}_{J_{s,t}}(\vphi)|
\prod_{\imath,\jmath}U^{\vphi\cdot
  W_{\imath,\jmath}}_{\imath,\jmath}\text{ and}\\
H_{w_i}(x)=&\sum_{y\in\T_{\al}\cap
  T_{\be}}y\sum_{\substack{\vphi\in\pi_2(x,y)\\\vphi\cdot
  Z=0\\\mu(\vphi)=0}}\vphi\cdot W_{i,s_i}|\mc{M}_{J_{s,t}}(\vphi)|U_{i,s_i}^{-1}
\prod_{\imath,\jmath}U^{\vphi\cdot
  W_{\imath,\jmath}}_{\imath,\jmath}.
\end{align*}

By analyzing the ends of $\mc{M}_{J_{s,t}}(\vphi)$ for Whitney disks
$\vphi$ with $\mu(\vphi)=1$, we see that
$F_{(\mc{H},J_s),(\mc{H}',J'_s)}\Psi_{i}+\Psi'_{i}F_{(\mc{H},J_s),(\mc{H}',J'_s)}=H_{z_i}\del+\del'
H_{z_i}$ and
$F_{(\mc{H},J_s),(\mc{H}',J'_s)}\Phi_{i}+\Phi'_{i}F_{(\mc{H},J_s),(\mc{H}',J'_s)}=H_{w_i}\del+\del'
H_{w_i}$.

\emph{\textbf{Case 2}: $J'_s=J_s$ and $\mc{H}'$ can be obtained from
$\mc{H}$ by isotoping and handlesliding
the $\al$ curves or by isotoping and handlesliding the $\be$ curves.} Without loss of generality, let us assume that we are isotoping and handlesliding the $\be$ curves. Furthermore, we can assume that the multicurves $\be$ and $\be'$ intersect each other transversely, $\wt{\mc{H}}=(\Si,\be,\be',z,w)$ is an admissible Heegaard diagram for the $n$-component unlink in $\#^g(S^1\times S^2)$, and the two tori $\T_{\be}$ and $\T_{\be'}$ intersect each other at $2^n$ points, all lying in the same Alexander grading. 

In $\wt{\mc{H}}$ let $\Theta\in\T_{\be}\cap\T_{\be'}$ be the element
with the highest Maslov grading. \cite[Lemma 9.1 and Lemma 9.4]{POZSz}
tell us that $\Theta$ is a cycle in $\CF_{(\wt{\mc{H}},J_s)}$. The
map $F_{(\mc{H},J_s),(\mc{H}',J_s)}$, evaluated on
$x\in\T_{\al}\cap\T_{\be}$, is given by
$$F_{(\mc{H},J_s),(\mc{H}',J_s)}(x)=\sum_{y\in\T_{\al}\cap\T_{\be'}}y\sum_{\substack{\vphi\in\pi_2(x,\Theta,y)\\\vphi\cdot
  Z=0\\\mu(\vphi)=0}}|\mc{M}_{J_s}(\vphi)|
\prod_{\imath,\jmath}U^{\vphi\cdot
  W_{\imath,\jmath}}_{\imath,\jmath},$$
where $\pi_2(x,\Theta,y)$ is the set of all Whitney triangles connecting $x$, $\Theta$ and $y$ \cite[Section 8.1.2]{POZSz}. Define $H_{z_i},H_{w_i}\colon\CF_{(\mc{H},J_s)}\rightarrow
\CF_{(\mc{H}',J_s)}$ as follows:
\begin{align*}
H_{z_i}(x)=&\sum_{y\in\T_{\al}\cap
  T_{\be'}}y\sum_{\substack{\vphi\in\pi_2(x,\Theta,y)\\\vphi\cdot
  Z=\vphi\cdot Z_i=1\\\mu(\vphi)=0}}|\mc{M}_{J_{s}}(\vphi)|
\prod_{\imath,\jmath}U^{\vphi\cdot
  W_{\imath,\jmath}}_{\imath,\jmath}\text{ and}\\
H_{w_i}(x)=&\sum_{y\in\T_{\al}\cap
  T_{\be'}}y\sum_{\substack{\vphi\in\pi_2(x,\Theta,y)\\\vphi\cdot
  Z=0\\\mu(\vphi)=0}}\vphi\cdot W_{i,s_i}|\mc{M}_{J_{s}}(\vphi)|U_{i,s_i}^{-1}
\prod_{\imath,\jmath}U^{\vphi\cdot
  W_{\imath,\jmath}}_{\imath,\jmath}.
\end{align*}

Since $\Theta$ is a cycle in $\CF_{(\wt{\mc{H}},J_s)}$ and all the
points in $\T_{\be}\cap\T_{\be'}$ lie in the same Alexander grading,
by counting the ends of $\mc{M}_{J_s}(\vphi)$ for Whitney triangles
$\vphi$ with $\mu(\vphi)=1$, we get
$F_{(\mc{H},J_s),(\mc{H}',J'_s)}\Psi_{i}+\Psi'_{i}F_{(\mc{H},J_s),(\mc{H}',J'_s)}=H_{z_i}\del+\del'
H_{z_i}$ and
$F_{(\mc{H},J_s),(\mc{H}',J'_s)}\Phi_{i}+\Phi'_{i}F_{(\mc{H},J_s),(\mc{H}',J'_s)}=H_{w_i}\del+\del'
H_{w_i}$.

\emph{\textbf{Case 3}: $\mc{H}'$ is obtained from $\mc{H}$ by an ordinary
(de)stabilization, as shown in \hyperref[fig:ordstab]{Figure \ref*{fig:ordstab}}, and $J'_s$ is related to $J_s$
as described below.} We fix a Riemann surface $E$ of genus $1$ with one $\al'$ circle and one $\be'$ circle, intersecting each other transversely at a single point. The Heegaard surface $\Si'$ is simply $\Si\# E$; $J'_s$ is induced from $J_s$, the complex structure on $E$, the two connected sum points in $\Si$ and $E$, and the length of the connected sum neck.

\begin{figure}
\psfrag{a}{$\al'$}
\psfrag{b}{$\be'$}
\psfrag{z}{}
\psfrag{p}{}
\begin{center}
\includegraphics[height=0.35\textwidth]{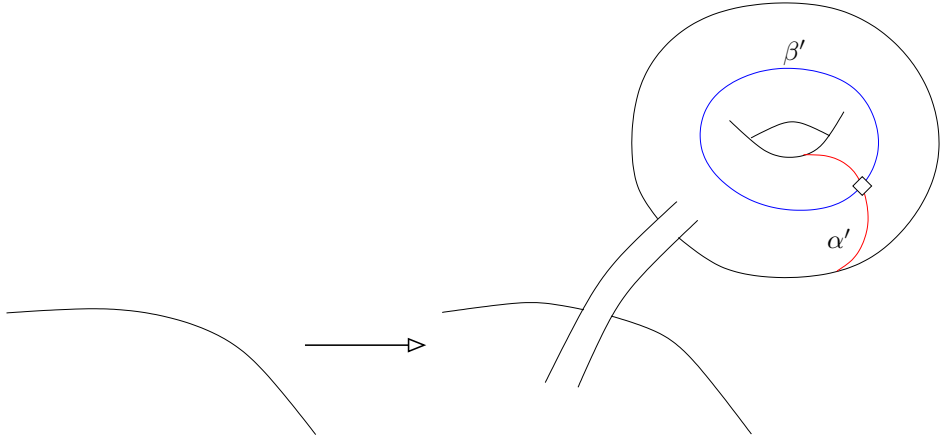}
\end{center}
\caption{An ordinary stabilization.}\label{fig:ordstab}
\end{figure}

There is a natural bijection between $\T_{\al}\cap\T_{\be}$ and $\T_{\al'}\cap\T_{\be'}$; let $x'\in\T_{\al'}\cap\T_{\be'}$ be the element corresponding to $x\in\T_{\al}\cap\T_{\be}$. For all $x,y\in\T_{\al}\cap\T_{\be}$, there is a natural bijection between $\pi_2(x,y)$ and $\pi_2(x',y')$ which preserves the Maslov index; let $\vphi'\in\pi_2(x',y')$ be the Whitney disk corresponding to $\vphi\in\pi_2(x,y)$.

However, by moving the connected sum points and by extending the
connected sum length, we can ensure that for all
$x,y\in\T_{\al}\cap\T_{\be}$ and for all $\vphi\in\pi_2(x,y)$ with
$\mu(\vphi)=1$, the two moduli spaces $\wh{\mc{M}}_{J_s}(\vphi)$ and
$\wh{\mc{M}}_{J'_s}(\vphi')$ are homeomorphic \cite[Theorem
10.4]{POZSz}. The map $F_{(\mc{H},J_s),(\mc{H}',J'_s)}$ sends $x$
to $x'$, and the map $F_{(\mc{H}',J'_s),(\mc{H},J_s)}$ is its
inverse. Therefore,
$F_{(\mc{H},J_s),(\mc{H}',J'_s)}\Psi_{i}+\Psi'_{i}F_{(\mc{H},J_s),(\mc{H}',J'_s)}
=F_{(\mc{H},J_s),(\mc{H}',J'_s)}\Phi_{i}+\Phi'_{i}F_{(\mc{H},J_s),(\mc{H}',J'_s)}=0$
and
$F_{(\mc{H}',J'_s),(\mc{H},J_s)}\Psi'_{i}+\Psi_{i}F_{(\mc{H}',J'_s),(\mc{H},J_s)}
=F_{(\mc{H}',J'_s),(\mc{H},J_s)}\Phi'_{i}+\Phi_{i}F_{(\mc{H}',J'_s),(\mc{H},J_s)}=0$.

\emph{\textbf{Case 4}: $\mc{H}'$ is obtained from $\mc{H}$ by a special
(de)stabilization, as shown in \hyperref[fig:spstab]{Figure \ref*{fig:spstab}}, and $J'_s$ is related to $J_s$
as described below.} We fix a Riemann surface $S$ of genus $0$ with one $\al'$ circle and one $\be'$ circle, intersecting each other transversely at two points $\varrho_1$ and $\varrho_2$. The Heegaard surface $\Si'$ is simply $\Si\# S$, where the connected sum is done near the $z_{\imath,\jmath}$ marking on $\Si$, and the $z$ and $w$ markings on $\Si'$ are as shown in \hyperref[fig:spstab]{Figure \ref*{fig:spstab}}; $J'_s$ is induced from $J_s$, the complex structure on $S$, the two connected sum points in $\Si$ and $S$, and the length of the connected sum neck. 

\begin{figure}
\psfrag{a}{$\al'$}
\psfrag{b}{$\be'$}
\psfrag{z}{$z_{\imath,\jmath}$}
\psfrag{z1}{$z_{\imath,k_{\imath}+1}$}
\psfrag{w}{$w_{\imath,k_{\imath}+1}$}
\psfrag{p}{$\varrho_1$}
\psfrag{q}{$\varrho_2$}
\begin{center}
\includegraphics[height=0.35\textwidth]{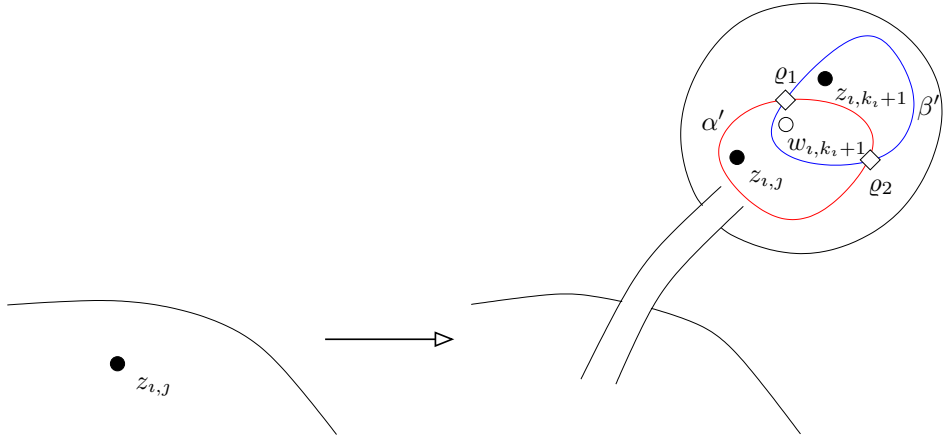}
\end{center}
\caption{A special stabilization.}\label{fig:spstab}
\end{figure}

There is a natural bijection between $\T_{\al'}\cap\T_{\be'}$ and
$(\T_{\al}\cap\T_{\be})\times\{1,2\}$, where the element corresponding
to $(x,1)$ uses the point $\varrho_1$ and is denoted by $x'$, and the
element corresponding to $(x,2)$ uses the point $\varrho_2$ and is
denoted by $x''$. Recall that $\CF_{(\mc{H}.J_s)}$ is the module
generated freely over $\mb{P}$ by $\T_{\al}\cap\T_{\be}$, and
$\CF_{(\mc{H}',J'_s)}$ is the module generated freely over
$\mb{P}'=\mb{P}\otimes\F_2[U_{\imath,k_{\imath}+1}]$ by
$\T_{\al'}\cap\T_{\be'}$. All the maps below are assumed to be
$\mb{P}$-module maps; furthermore, the maps between $\mb{P}'$-modules
are assumed to be $\mb{P}'$-module maps. \cite[Proposition
6.5]{POZSzlinkinvariants} tells us that by moving the connected sum
point in $S$ near the $\al'$ circle and by extending the connected sum
length, we can ensure the following: $\CF_{(\mc{H}',J'_s)}$ is
isomorphic to the mapping
cone 
\xymatrixcolsep{5pc}\xymatrix{C_2\otimes\F_2[U_{\imath,k_{\imath}+1}]\ar[r]^{U_{\imath,k_{\imath}+1}+U_{\imath,a_{\imath}(\jmath)}}
&C_1 \otimes\F_2[U_{\imath,k_{\imath}+1}]}, where $C_1$ is isomorphic
to $\CF_{(\mc{H},J_s)}$ corresponding to elements of the form $x'$,
and $C_2$ is isomorphic to $\CF_{(\mc{H},J_s)}$ corresponding to
elements of the form $x''$;
$F_{(\mc{H},J_s),(\mc{H}',J'_s)}(x)=x'$;
$F_{(\mc{H}',J'_s),(\mc{H},J_s)}(U_{\imath,k_{\imath}+1}^m x'')=0$
and $F_{(\mc{H}',J'_s),(\mc{H},J_s)}(U_{\imath,k_{\imath}+1}^m
x')=U_{\imath,a_{\imath}(\jmath)}^m x$; for
$(i,j)\notin\{(\imath,\jmath),(\imath,k_{\imath}+1)\}$,
$\Psi'_{{i,j}}(x')=(\Psi_{{i,j}}(x))'$ and
$\Psi'_{{i,j}}(x'')=(\Psi_{{i,j}}(x))''$;
$(\Psi'_{{\imath,\jmath}}+\Psi'_{{\imath,k_{\imath}+1}})(x')=(\Psi_{{\imath,\jmath}}(x))'$
and
$(\Psi'_{{\imath,\jmath}}+\Psi'_{{\imath,k_{\imath}+1}})(x'')=(\Psi_{{\imath,\jmath}}(x))''$;
$\Phi'_{{i}}(x')=(\Phi_{{i}}(x))'$ and
$\Phi'_{{i}}(x'')=(\Phi_{{i}}(x))''+x'\delta_{i\imath}\delta_{s_i
a_{\imath}(\jmath)}$.

Therefore,
$F_{(\mc{H},J_s),(\mc{H}',J'_s)}\Psi_{i}+\Psi'_{i}F_{(\mc{H},J_s),(\mc{H}',J'_s)}
=F_{(\mc{H},J_s),(\mc{H}',J'_s)}\Phi_{i}+\Phi'_{i}F_{(\mc{H},J_s),(\mc{H}',J'_s)}=0$. Similarly,
we get
$F_{(\mc{H}',J'_s),(\mc{H},J_s)}\Psi'_{i}+\Psi_{i}F_{(\mc{H}',J'_s),(\mc{H},J_s)}=0$
and
$(F_{(\mc{H}',J'_s),(\mc{H},J_s)}\Phi'_{i}+\Phi_{i}F_{(\mc{H}',J'_s),(\mc{H},J_s)})(U_{\imath,k_{\imath}+1}^m
x')=0$. However,
$(F_{(\mc{H}',J'_s),(\mc{H},J_s)}\Phi'_{i}+\Phi_{i}F_{(\mc{H}',J'_s),(\mc{H},J_s)})(U_{\imath,k_{\imath}+1}^m
x'')=U_{\imath,a_{\imath}(\jmath)}^mx\delta_{i\imath}\delta_{s_i
a_{\imath}(\jmath)}$. Define
$H_{w_i}\colon\CF_{(\mc{H}',J'_s)}\rightarrow \CF_{(\mc{H},J_s)}$ as
follows: $H_{w_i}(U_{\imath,k_{\imath}+1}^m x'')=0$ and
$H_{w_i}(U_{\imath,k_{\imath}+1}^m x')=m
U_{\imath,a_{\imath}(\jmath)}^{m-1} x\delta_{i\imath}\delta_{s_i
a_{\imath}(\jmath)}$. A careful analysis shows that
$F_{(\mc{H}',J'_s),(\mc{H},J_s)}\Phi'_{i}+\Phi_{i}F_{(\mc{H}',J'_s),(\mc{H},J_s)}=H_{w_i}\del'+\del
H_{w_i}$. Observe that the chain homotopy is still $U_i$-equivariant.
\end{proof}

We will now state and prove some properties of the maps $\Psi_{i}$ and $\Phi_{i}$.

\begin{lem}\label{lem:homotopyproperties}
  For all $i$, all three of the maps
  $\Psi_{i}\Phi_{i}+\Phi_{i}\Psi_{i}$, $\Psi^2_{i}$ and $\Phi^2_i$ are
  chain homotopic to zero, where the chain homotopies are also
  $U_{i,j}$-equivariant.
\end{lem}

\begin{proof}
Define the three $U_{i,j}$-equivariant chain homotopies as follows:
\begin{align*}
  H_{1,i}(x)=&\sum_{y\in\T_{\al}\cap
    T_{\be}}y\sum_{\substack{\vphi\in\pi_2(x,y)\\\vphi\cdot
      Z=\vphi\cdot Z_i=1\\\mu(\vphi)=1}}(\vphi\cdot
  W_{i,s_i})|\wh{\mc{M}}_{J_s}(\vphi)|U^{-1}_{i,s_i}
  \prod_{\imath,\jmath}U^{\vphi\cdot
    W_{\imath,\jmath}}_{\imath,\jmath},\\
  H_{2,i}(x)=&\sum_{y\in\T_{\al}\cap
    T_{\be}}y\sum_{\substack{\vphi\in\pi_2(x,y)\\\vphi\cdot
      Z=\vphi\cdot Z_i=2\\\mu(\vphi)=1}}|\wh{\mc{M}}_{J_s}(\vphi)|
  \prod_{\imath,\jmath}U^{\vphi\cdot
    W_{\imath,\jmath}}_{\imath,\jmath}\text{ and}\\
  H_{3,i}(x)=&\sum_{y\in\T_{\al}\cap\T_{\be}}y
  \sum_{\substack{\vphi\in\pi_2(x,y)\\\vphi\cdot Z=0\\\mu(\vphi)=1}}
  {{\vphi\cdot W_{i,s_i}}\choose 2}|\wh{\mc{M}}_{J_s}(\vphi)|U^{-2}_{i,s_i}
  \prod_{\imath,\jmath}U^{\vphi\cdot
    W_{\imath,\jmath}}_{\imath,\jmath} .
\end{align*}

By counting the ends of $\wh{\mc{M}}_{J_s}(\vphi)$ for Whitney disks
$\vphi$ with $\mu(\vphi)=2$ and $\vphi\cdot Z=\vphi\cdot Z_i=1$, we
see that $\Psi_{i}\Phi_{i}+\Phi_{i}\Psi_{i}=\del
H_{1,i}+H_{1,i}\del$. When $\vphi\cdot W_i=1$, we might have boundary
degenerations in the ends of $\wh{\mc{M}}_{J_s}(\vphi)$; however, they
cancel in pairs \cite[Theorem 5.5]{POZSzlinkinvariants}. Similarly, by
counting the ends of $\wh{\mc{M}}_{J_s}(\vphi)$ for Whitney disks
$\vphi$ with $\mu(\vphi)=2$ and $\vphi\cdot Z=\vphi\cdot Z_i=2$, we
see that $\Psi^2_{i}=\del H_{2,i}+H_{2,i}\del$. For the third case,
let us explicitly do the calculation.
\begin{align*}
&\qquad([H_{3,i}\colon\del]+\Phi_i^2)(x)\\
=&\sum_{y\in\T_{\al}\cap
  T_{\be}}y\sum_{\substack{\vphi\in\pi_2(x,y)\\\vphi\cdot
  Z=0\\\mu(\vphi)=2}}
\sum_{\substack{\vphi_1,\vphi_2,m_1,m_2\\\vphi_1\cdot W_{i,s_i}=m_1\\
    \vphi_2\cdot W_{i,s_i}=m_2\\\mu(\vphi_1)=\mu(\vphi_2)=1\\ \vphi=\vphi_1*\vphi_2}}
({m_1\choose 2}+{m_2\choose 2}+m_1m_2)|\wh{\mc{M}}_{J_s}(\vphi_1)\times\wh{\mc{M}}_{J_s}(\vphi_2)|
  U^{-2}_{i,s_i} \prod_{\imath,\jmath}U^{\vphi\cdot
    W_{\imath,\jmath}}_{\imath,\jmath}\\
=&\sum_{y\in\T_{\al}\cap
  T_{\be}}y\sum_{\substack{\vphi\in\pi_2(x,y)\\\vphi\cdot
  Z=0\\\mu(\vphi)=2}}
\sum_{\substack{\vphi_1,\vphi_2,m_1,m_2\\\vphi_1\cdot W_{i,s_i}=m_1\\
    \vphi_2\cdot W_{i,s_i}=m_2\\\mu(\vphi_1)=\mu(\vphi_2)=1\\ \vphi=\vphi_1*\vphi_2}}
{{m_1+m_2}\choose 2}|\wh{\mc{M}}_{J_s}(\vphi_1)\times\wh{\mc{M}}_{J_s}(\vphi_2)|
  U^{-2}_{i,s_i} \prod_{\imath,\jmath}U^{\vphi\cdot W_{\imath,\jmath}}_{\imath,\jmath}\\
=&\sum_{y\in\T_{\al}\cap
  T_{\be}}y\sum_{\substack{\vphi\in\pi_2(x,y)\\\vphi\cdot
  Z=0\\\mu(\vphi)=2}}{{\vphi\cdot W_{i,s_i}}\choose 2}U_{i,s_i}^{-2}
\sum_{\substack{\vphi_1,\vphi_2\\\mu(\vphi_1)=\mu(\vphi_2)=1\\ \vphi=\vphi_1*\vphi_2}}
|\wh{\mc{M}}_{J_s}(\vphi_1)\times\wh{\mc{M}}_{J_s}(\vphi_2)|
  \prod_{\imath,\jmath}U^{\vphi\cdot W_{\imath,\jmath}}_{\imath,\jmath}\\
=&0\qedhere.
\end{align*}
\end{proof}

\begin{thm}
For all $i$, the map $\Id+\Psi_{i}\Phi_{i}=\Id+\Phi_i\Psi_i$ is a link-invariant
involution of $\CFL(Y,L,p)$, viewed as an object in $N(K(\mc{A}_l))$.
\end{thm}

\begin{proof}
We already know
from \hyperref[thm:invariance]{Theorem \ref*{thm:invariance}} that
$\Id+\Psi_{i}\Phi_{i}$ is a link-invariant map from
$\CFL(Y,L,p)$ to itself. To see that it is an involution, observe
that in $N(K(\mc{A}_l))$, as a consequence
of \hyperref[lem:homotopyproperties]{Lemma \ref*{lem:homotopyproperties}},
$(\Id+\Psi_{i}\Phi_{i})^2=\Id+\Psi_{i}\Phi_{i}\Psi_{i}\Phi_{i}=
\Id+\Psi_{i}\Psi_{i}\Phi_{i}\Phi_{i}=\Id$.
\end{proof}

Instead of working in the slightly unfamiliar category
$N(K(\mc{A}_l))$, we often take the homology and work with
$\HFL(Y,L,p)=H_*(\CFL(Y,L,p))$ in $N(\mc{B}_l)$. Another standard
object to work with is the hat invariant $\wh{\HFL}(Y,L,p)$ living in
$N(\mc{C}_l)$: it is the homology of the mapping cone of all the maps
$U_1,\ldots,U_l$ in $\CFL(Y,L,p)$.

The map $\Id+\Psi_{i}\Phi_{i}$ induces link-invariant involutions on
$\HFL(Y,L,p)$ and $\wh{\HFL}(Y,L,p)$. Even though the involution on
$\wh{\HFL}(Y,L,p)$ (and hence the one on $\CFL(Y,L,p)$) is often
non-trivial, cf.\ \hyperref[thm:computations]{Theorem
  \ref*{thm:computations}}, quite (un-)surprisingly, the involution on
$\HFL(Y,L,p)$ is always the identity.

\begin{lem}\label{lem:triviality}
For all $i$, the involution $\Id+\Psi_{i}\Phi_{i}$ on
$\HFL(Y,L,p)$, viewed as an object in $N(\mc{B}_l)$, is
the identity map.
\end{lem}

\begin{proof}
In order to prove this, we only need to show that the map
$\Psi_{i}\Phi_{i}$ is chain homotopic to zero, where the chain
homotopy need not be $U_i$-equivariant. Fix a Heegaard diagram
$\mc{H}$ and a path of nearly symmetric almost complex structures
$J_s$ on the symmetric product. The chain homotopy $H_i$ is
$U_{\imath,\jmath}$-equivariant for all $(\imath,\jmath)\neq (i,s_i)$,
and is defined as follows:
$$H_i(U^m_{i,s_i}x)=\sum_{y\in\T_{\al}\cap\T_{\be}}y\sum_{\substack{\vphi\in\pi_2(x,y)\\\vphi\cdot
  Z=\vphi\cdot Z_i=1\\\mu(\vphi)=1}}m|\wh{\mc{M}}_{J_s}(\vphi)|U^{m-1}_{i,s_i}
\prod_{\imath,\jmath}U^{\vphi\cdot
  W_{\imath,\jmath}}_{\imath,\jmath}.$$

A careful analysis of the ends of $\wh{\mc{M}}_{J_s}(\vphi)$ for
Whitney disks $\vphi$ with $\mu(\vphi)=2$ and $\vphi\cdot Z=\vphi\cdot
Z_i=1$ shows that $\Psi_{i}\Phi_{i}=\del H_i+H_i\del$.
\end{proof}

In view of \hyperref[lem:triviality]{Lemma \ref*{lem:triviality}}
along with \hyperref[thm:main]{Theorem \ref*{thm:main}}, we can define
the invariant $\HFL$ for unpointed links in $S^3$, although due to
\hyperref[thm:computations]{Theorem \ref*{thm:computations}}, we must
continue to treat $\CFL$ and $\wh{\HFL}$ as invariants of pointed
links. Let us conclude this section with a rather bold conjecture,
which is somewhat justified by \hyperref[thm:main]{Theorem
  \ref*{thm:main}}.

\begin{conj}
The two automorphisms $\rho(\si_i)$ and $\Id+\Psi_{i}\Phi_{i}$ in
$\Aut_{N(K(\mc{A}_l))}(\CFL(Y,L,p))$ are equal.
\end{conj}

\section{Grid diagrams}\label{sec:grids}
In this section, we concentrate on pointed links in $S^3$. The main
aim is to prove \hyperref[thm:main]{Theorem \ref*{thm:main}}. We focus
our attention on the component $(L_1,p_1)$ of the $l$-component
pointed link $(L,p)$ in $S^3$. Following \cite{CMPOSS}, we will
represent such links by a special type of Heegaard diagrams
called \emph{grid diagrams}. A grid diagram of index $n$ is a Heegaard
diagram $\mc{G}=(T,\al,\be,z,w)$ for $(L,p)$, where the Heegaard
surface $T$ is the torus obtained as a quotient of $[0,1]^2\sbs \C$ by
identifying opposite sides, $\al$ is a multicurve which is isotopic to
the image of $[0,1]\times \{0,\frac{1}{n},\ldots,\frac{n-1}{n}\}$,
$\be$ is a multicurve which is isotopic to the image of
$\{0,\frac{1}{n},\ldots,\frac{n-1}{n}\}\times [0,1]$, and each
$\al$-circle intersects each $\beta$-circle at exactly one point. It
is easy to see that each pointed link can be represented by a grid
diagram. We usually take $J_s$ to be the constant path of the product
complex structure on $\Sym^n(T)$ induced from the complex structure on
$\C$. By generically perturbing $\al$ and $\be$, we can ensure that
$J_s$ achieves transversality \cite[Proposition 3.9]{RL}. We will keep
using our notations from the previous sections. However, to avoid
clutter, from now on, unless we deem it to be particularly
illuminating, we will drop the subscript $J_s$ from our notation. We
sometimes use the words north, south, east and west to denote local
directions on the torus $T$ (i.e.\ directions on some contractible
subset of $T$); at all such times, it is implicitly understood that we
have isotoped the $\al$ and $\be$ circles to horizontal and vertical
circles respectively.

We will set up for the proof, while doing a brief review of grid diagrams, in the following few subsections. 

\subsection{The grid chain complex}
Given a grid diagram $\mc{G}$ of index $n$, there are exactly $n!$ points in $\T_{\al}\cap\T_{\be}$. If $\vphi\in\pi_2(x,y)$ for some $x,y\in\T_{\al}\cap\T_{\be}$, then let $D(\vphi)$ be the shadow of $\vphi$ which is a $2$-chain generated by components of $T\sm(\al\cup\be)$ \cite[Definition 2.13]{POZSz}. The domain $D(\vphi)$ is said to be non-negative if it has non-negative coefficients everywhere. Since we are working with a product complex structure, if $\mc{M}_{J_s}(\vphi)\neq\varnothing$, then $D(\vphi)$ is non-negative. The central fact about grid diagrams is the following observation.

\begin{thm}\cite{CMPOSS}\label{thm:cmpossmain}
If $\vphi\in\pi_2(x,y)$ is a Whitney disk in a grid diagram with $D(\vphi)$ non-negative, then $\mu(\vphi)\geq 0$. Furthermore, $\mu(\vphi)=0$ happens precisely when $x=y$ and $\vphi$ is the trivial disk; and $\mu(\vphi)=1$ happens precisely when $x$ and $y$ differ in exactly two coordinates, $D(\vphi)$ is a properly embedded rectangle in $(T,\al\cup\be)$ which does not contain any coordinates of $x$ or $y$ in its interior, the northeast and southwest corners of $D(\vphi)$ are coordinates of $x$ and the northwest and southeast corners of $D(\vphi)$ are coordinates of $y$, and in that case $\wh{\mc{M}}(\vphi)$ has exactly one point.
\end{thm}

In particular, this theorem implies that the grid chain complex $\CF_{(\mc{G},J_s)}\in \Ob_{K(\mc{A}_l)}$ is independent of $J_s$, as long as it is a constant path of the product complex structure induced from some complex structure on $T$. In fact, the following is an explicit description of the chain complex $\CF_{\mc{G}}$ in grid diagram terminology.

A \emph{state} $x$ is an $n$-tuple of points $x=(x_1,\ldots,x_n)$ (and the points $x_i$ are called the \emph{coordinates} of $x$), such that each $\al$-circle contains some $x_i$ and each $\be$ circle contains some $x_j$. Clearly there are $n!$ states, and there is a natural bijection between $\T_{\al}\cap\T_{\be}$ and the set of all states $\mc{S}_{\mc{G}}$. A \emph{grid $2$-chain} is a formal linear combination of the $n^2$ components of $T\sm(\al\cup\be)$ over $\Z$. A grid $2$-chain is said to be \emph{positive} if all its coefficients are non-negative. Given a point $p\in T\sm(\al\cup\be)$ and a grid $2$-chain $D$, the number $n_p(D)$ is the coefficient of $D$ at the component of $T\sm(\al\cup\be)$ that contains the point $p$. For any grid $2$-chain $D$, let $n_{z_i}(D)=\sum_j n_{z_{i,j}}(D)$, $n_{w_i}(D)=\sum_j n_{w_{i,j}}(D)$, $n_{z}(D)=\sum_i n_{z_{i}}(D)$ and $n_{w}(D)=\sum_i n_{w_{i}}(D)$. A \emph{domain} joining a state $x$ to a state $y$ is a grid $2$-chain $D$ such that $\del(\del D\cap\al)=y-x$. The set of all domains joining $x$ to $y$ is denoted by $\mc{D}_{\mc{G}}(x,y)$, and it is in natural bijection with $\pi_2(x,y)$. A \emph{rectangle} $R\in\mc{D}_{\mc{G}}(x,y)$ is a domain satisfying the following conditions: each coefficient of $R$ is either $0$ or $1$; the closure of the region where $R$ has coefficient $1$ is properly embedded rectangle in $T$; that rectangle does not contain any coordinates of $x$ or $y$ in its interior; the northeast and southwest corners of that rectangle are coordinates of $x$ and the northwest and southeast corners of that rectangle are coordinates of $y$. The set of all rectangles joining $x$ to $y$ is denoted by $\mc{R}_{\mc{G}}(x,y)\sbs \mc{D}_{\mc{G}}(x,y)$.

The $(M,A_1,\ldots,A_l)$-graded ring $\mb{P}$ is the polynomial ring generated over $\F_2$ by the variables $U_{i,j}$ for $i\in\{1,\ldots,l\},j\in\{1,\ldots,k_i\}$, and $\CF_{\mc{G}}$ is the $\F_2[U_1,\ldots,U_l]$-module freely generated over $\mb{P}$ by $\mc{S}_{\mc{G}}$, where the $U_i$-action is multiplication by $U_{i,s_i}$. The $U_{i,j}$-equivariant $(-1,0,\ldots,0)$-graded boundary map, evaluated on a state $x$, is
$$\del x=\sum_{y\in\mc{S}_{\mc{G}}}y \sum_{\substack{R\in\mc{R}_{\mc{G}}(x,y)\\ n_z(R)=0}}\prod_{\imath,\jmath} U_{\imath,\jmath}^{n_{w_{\imath,\jmath}}(R)}.$$

\subsection{Changing the complex structure} Let $J_s$ and $J'_s$ be
the constant paths of almost complex structures on $\Sym^n(T)$ induced
from two complex structure on $T$. We know that
$\CF_{(\mc{G},J_s)}$ and $\CF_{(\mc{G},J'_s)}$ are
the \emph{same} object; therefore, it is not unnatural to expect the
naturality map
$F_{(\mc{G},J_s),(\mc{G},J'_s)}\in\Mor_{K(\mc{A}_l)}(\CF_{(\mc{G},J_s)},\CF_{(\mc{G},J'_s)})$
to be the identity map. This is indeed the case.

\begin{thm}
If $J_s$ and $J'_s$ are two constant paths of almost complex
structures on $\Sym^n(T)$ induced from two complex structures on $T$,
then $F_{(\mc{G},J_s),(\mc{G},J'_s)}=\Id.$
\end{thm}

\begin{proof}[Sketch of a proof]
This is a direct consequence of \cite[Proof of Theorem 6.1]{POZSz},
adapted to our present setting, where any Whitney disk
$\vphi\in\pi_2(x,y)$ whose shadow $D(\vphi)\in\mc{D}_{\mc{G}}(x,y)$ is
positive and which satisfies $\mu(\vphi)\leq 0$, is a trivial one.
\end{proof}

\subsection{Commutation} Commutation comes in two flavors, horizontal commutation and vertical commutation. The story for vertical commutation can be guessed from the story of horizontal commutation by reversing the roles of $\al$ and $\be$; vertical commutations are also described in full detail in \cite{CMPOZSzDT}, from where much of the material for this subsection is derived; so for now, let us only talk about horizontal commutation.

\begin{figure}
\psfrag{a}{$\al$}
\psfrag{p}{$\varrho$}
\psfrag{p'}{$\varrho'$}
\psfrag{a1}{$\al_1$}
\psfrag{a11}{$\al'_1$}
\begin{center}
\includegraphics[width=0.5\textwidth]{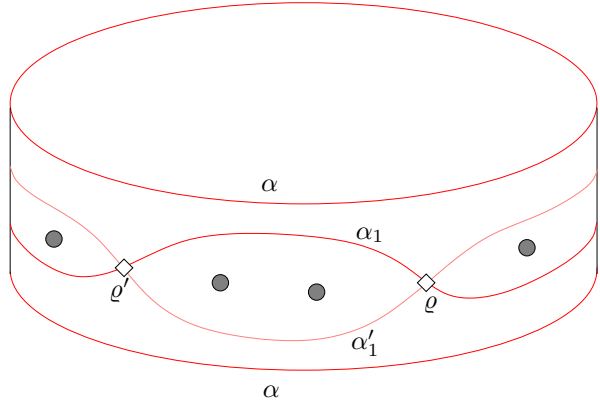}
\end{center}
\caption{Horizontal commutation.}\label{fig:commutation}
\end{figure}

A horizontal commutation is a pair of grid diagrams $(\mc{G},\mc{G}')$ drawn on the same torus, such that $\mc{G}'$ can be obtained from $\mc{G}$ by changing exactly one $\al$-circle in the following manner: if the circle $\al_1$ in $\mc{G}$ is changed to the circle $\al'_1$ in $\mc{G}'$, then some neighborhood of $\al_1\cup\al'_1$ must be homeomorphic to the region shown in \hyperref[fig:commutation]{Figure \ref*{fig:commutation}}, where the four black dots are two $z$ markings and two $w$ markings. The two points in $\al_1\cap\al'_1$ are marked as $\varrho$ and $\varrho'$. The $\be$-circles are not shown; but we assume that they avoid both $\varrho$ and $\varrho'$.

Commutation induces the naturality map
$F_{\mc{G},\mc{G}'}\in\Mor_{K(\mc{A}_l)}(\CF_{\mc{G}},\CF_{\mc{G}'})$.
Allow us to explain this map in
grid diagram terminology. 

For $x\in\mc{S}_{\mc{G}}$ and
$y\in\mc{S}_{\mc{G}'}$, a \emph{pentagon} $P$ joining $x$ to $y$ is a
$2$-chain generated by the components of $T\sm(\al\cup\al'\cup\be)$,
such that: each component appears with coefficient $0$ or $1$; the
closure of the union of the coefficient $1$ components is an embedded
pentagon whose one of the vertices is $\varrho$; the embedded pentagon
does not contain any coordinates of $x$ or $y$ in its interior; the
northeast and southwest corners of the embedded pentagon are
coordinates of $x$; and the northwest and southeast corners of the
embedded pentagon are coordinates of $y$. The set of all pentagons
joining $x$ to $y$ is denoted by
$\mbox{${}_{\varrho}\mc{P}_{\mc{G},\mc{G}'}(x,y)$}$. We have an
$U_{i,j}$-equivariant $(0,0,\ldots,0)$-graded chain map
$\mbox{${}_{\varrho}\wt{F}_{\mc{G},\mc{G}'}$} \colon \CF_{\mc{G}} \rightarrow \CF_{\mc{G}'}$,
which when evaluated on $x\in\mc{S}_{\mc{G}}$ is given by
$$\mbox{${}_{\varrho}\wt{F}_{\mc{G},\mc{G}'}(x)$}=
\sum_{y\in\mc{S}_{\mc{G}'}}y\sum_{\substack{P\in{}_{\varrho}\mc{P}_{\mc{G},\mc{G}'}(x,y)\\ n_z(P)=0}}\prod_{\imath,\jmath} U_{\imath,\jmath}^{n_{w_{\imath,\jmath}}(P)}.$$

This chain map depends on the location of the point $\varrho$, namely
which component of $T\sm\be$ contains $\varrho$. Similarly, we get
another chain map
$\mbox{${}_{\varrho'}\wt{F}_{\mc{G}',\mc{G}}$} \colon \CF_{\mc{G}'} \rightarrow \CF_{\mc{G}}$,
which depends on the location of the point
$\varrho'$. \cite[Proposition 3.2]{CMPOZSzDT} tells us that any of the maps
$\mbox{${}_{\varrho'}\wt{F}_{\mc{G}',\mc{G}}$}$ is an inverse for any
of the maps $\mbox{${}_{\varrho}\wt{F}_{\mc{G},\mc{G}'}$}$ in the
homotopy category $K(\mc{A}_l)$. Therefore, we get a well-defined map
in $\wt{F}_{\mc{G},\mc{G}'}\in \Mor_{K(\mc{A}_l)}(\CF_{\mc{G}},\CF_{\mc{G}'})$
which is independent of $\varrho$. 

It is not hard to check that this
map $\wt{F}_{\mc{G},\mc{G}'}$, defined in terms of pentagons in the torus $T$, is same as the map
$F_{(\mc{G},J_s),(\mc{G}',J_s)}$, defined using holomorphic
triangles in $\Sym^n(T)$. Anyone who believes in the truth of the above assertion should fearlessly skip the rest of this subsection.

\begin{figure}
\psfrag{a}{$\al$}
\psfrag{b}{$\be$}
\psfrag{a'}{$\al'$}
\begin{center}
\includegraphics[width=0.6\textwidth]{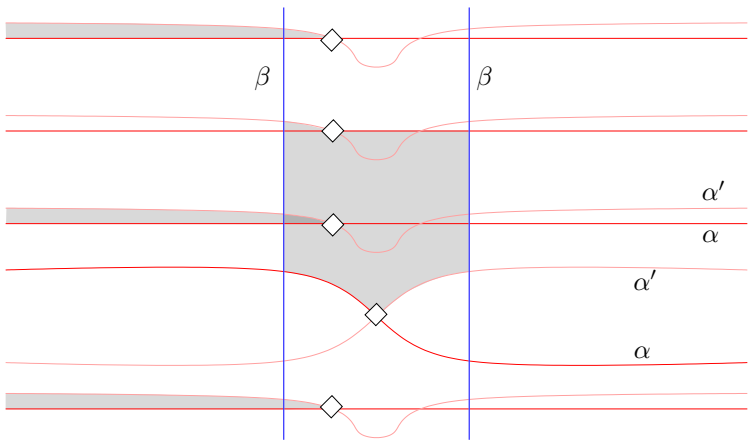}
\end{center}
\caption{Horizontal commutation.}\label{fig:pentagonholo}
\end{figure}

Given a commutation diagram, we can perturb $\al'$ to ensure that each $\al'$-circle intersects its corresponding $\al$-circle in exactly two points; furthermore, except for $\varrho'$, we can ensure that each of the intersection points lies in the component of $T\sm\be$ that contains the point $\varrho$. The situation is illustrated in \hyperref[fig:pentagonholo]{Figure \ref*{fig:pentagonholo}}. Let $\Theta$ be the top dimensional generator in the Heegaard diagram $(T,\al',\al)$; the coordinates of $\Theta$ are shown (note, $\varrho$ is one of the coordinates).

For $x\in\mc{S}_{\mc{G}}$ and $y\in\mc{S}_{\mc{G}'}$, if  $D$ is a $2$-chain generated by the components of $T\sm(\al\cup\al'\cup\be)$ such that $\del(\del D|_{\al})=\Theta-x$ and $\del(\del D|_{\al'})=y-\Theta$, then $D$ is called a domain joining $x$, $\Theta$ and $y$. If $\vphi\in\pi_2(x,\Theta,y)$ is a Whitney triangle connecting $x$, $\Theta$ and $y$ in $\Sym^n(T)$, then its shadow $D(\vphi)$ is a domain joining $x$, $\Theta$ and $y$. Conversely, given any domain $D$ joining $x$, $\Theta$ and $y$, there is a unique Whitney triangle $\vphi\in\pi_2(x,\Theta,y)$ such that $D=D(\vphi)$.

It is easy to see that each pentagon in
$\mbox{${}_{\varrho}\mc{P}_{\mc{G},\mc{G}'}(x,y)$}$ gives rise to
exactly two positive domains joining $x$, $\Theta$ and $y$. Let us
call such domains twin pentagonal domains. One such pentagonal domain is
shown in \hyperref[fig:pentagonholo]{Figure \ref*{fig:pentagonholo}}.

\begin{lem}
If $\vphi\in\pi_2(x,\Theta,y)$ is a Whitney triangle with
$\mu(\vphi)=0$ and $D(\vphi)$ positive, then $D(\vphi)$ is a
pentagonal domain. Conversely, if $\vphi,\vphi'\in\pi_2(x,\Theta,y)$
are Whitney triangles such that $D(\vphi)$ and $D(\vphi')$ are twin
pentagonal domains, then $\mu(\vphi)=\mu(\vphi')=0$ and
$|\mc{M}(\vphi)|+|\mc{M}(\vphi')|=1$. Therefore, the triangle map
$F_{\mc{G},\mc{G}'}$ agrees with $\wt{F}_{\mc{G},\mc{G}'}$.
\end{lem}

\begin{proof}
For the first direction, assume $\vphi\in\pi_2(x,\Theta,y)$ is a
Whitney triangle with $\mu(\vphi)=0$ and $D(\vphi)$ positive. Let $B$
be the union of the $2n$ bigonal components of
$T\sm(\al\cup\al')$. Let $D'$ be the unique domain joining some state
$y'\in\mc{S}_{\mc{G}}$ and $\Theta$ and $y$, such that $D'$ is
supported in $B$, and $\del D'|_{\al'}=\del D(\vphi)|_{\al'}$. Using
the Maslov index formula from \cite{SS}, it is easy to check that the
Maslov index of $D'$ is $(-1)$. Therefore, $\wt{D}=D(\vphi)-D'$ is a
Maslov index $1$ domain in $\mc{D}_{\mc{G}}(x,y')$. Furthermore, since
$D'$ has coefficient zero outside $B$, and $D(\vphi)$ is positive,
$\wt{D}$ is a positive domain as
well. \hyperref[thm:cmpossmain]{Theorem \ref*{thm:cmpossmain}} tells
us that $\wt{D}$ is a rectangle in $\mc{G}$. However, we know that
$D(\vphi)=\wt{D}+D'$ is also a positive domain. Therefore, the
rectangle $\wt{D}$ can only be in certain configurations. By analyzing
them carefully, we see that those are precisely the configurations for
which $D(\vphi)$ is a pentagonal domain.

The other direction is fairly straightforward. Let
$\vphi,\vphi'\in\pi_2(x,\Theta,y)$ be Whitney triangles such that
$D(\vphi)$ and $D(\vphi')$ are twin pentagonal domains. A direct
computation reveals that the Maslov index of any pentagonal domain is
zero. To show that $|\mc{M}(\vphi)|+|\mc{M}(\vphi')|=1$, following
standard practice, we will use Lipshitz's cylindrical
reformulation \cite{RL}; in the cylindrical version, we will count the
number of holomorphic embeddings of surfaces $F\hookrightarrow
T\times \{z\in\C\mid |z|\leq 1,z^3\neq 1\}$ satisfying certain
conditions, such as: the image of the projection onto the first factor
is either $D(\vphi)$ or $D(\vphi')$, and the projection onto the
second factor is an $n$-sheeted cover of $\{z\in\C\mid |z|\leq
1,z^3\neq 1\}$ with exactly one branch point (the number of branch
points can be figured out from another formula
in \cite{SS}). Therefore, $F$ must be a disjoint union of $(n-2)$
copies of a disk punctured at three boundary points (call them
$3$-gons), and a single copy of a disk punctured at six boundary
points (call it a $6$-gon). Therefore, the moduli space
$\mc{M}(\vphi)\sqcup\mc{M}(\vphi')$ is the product of the $(n-1)$ moduli spaces coming
from these $(n-1)$ components. It is a fairly easy exercise in complex
analysis to check that the moduli space of embeddings of a $3$-gon
contains exactly one point.

To show that the moduli space of embeddings of the $6$-gon into
$T\times \{z\in\C\mid |z|\leq 1,z^3\neq 1\}$ contains an odd number of
points, let us do a model calculation. If $\varrho$ lies in the
southern portion of the pentagonal domains $D(\vphi)$ and
$D(\vphi')$, then one of $D(\vphi)$ and $D(\vphi')$ looks like the
shaded hexagonal region
in \hyperref[fig:pentagonholo]{Figure \ref*{fig:pentagonholo}}, while
the other one is its twin; if $\varrho$ lies in the northern portion
of the pentagonal domain, then $D(\vphi)$ and $D(\vphi')$ have similar
but different shapes. Call $D(\vphi)$ and $D(\vphi')$ the original
twin pentagonal domains. Now consider the index $2$ grid diagram $\mc{G}_0$
for the $2$-component unlink. There are two states and they lie in
different Maslov gradings. Do a horizontal commutation to obtain the
grid diagram $\mc{G}_1$, which also has two states, lying in different
Maslov gradings. The triangle map must be an isomorphism, and we have
already seen that the shadow of each Whitney triangle must be a
pentagonal domain. There are exactly two sets of twin pentagonal
domains in this model commutation diagram. Consider the twin
pentagonal domains that have
the same shape as the original twin pentagonal domains. Call them the
model twin pentagonal domains. Choose a complex structure on the model
twin pentagonal domains which matches the one on the original twin pentagonal
domains. Since the triangle map for the model commutation is a graded
isomorphism, the moduli space corresponding to the model twin pentagonal
domains, and hence the moduli space corresponding to the original
twin pentagonal domains, must contain an odd number of points.
\end{proof}

\subsection{Stabilization and destabilization}
Stabilization for grid diagrams is a move that converts a grid diagram
of index $n$ to a grid diagram of index $(n+1)$. There are several
variants (four or eight, depending on how we count them) of
stabilization. However, we are mostly concerned with only one of these
configurations, so let us describe it in detail.

Let $\mc{G}=(T,\al,\be,z,w)$ be an index $n$ grid diagram where the
link component $L_1$ contains the $2k_1$ markings
$z_{1,1},\ldots,\allowbreak
z_{1,k_1},w_{1,1},\allowbreak \ldots,w_{1,k_1}$. Let
$\mc{G}'=(T,\al',\be',z',w')$ an index $(n+1)$ grid diagram
representing $L$ and satisfying the following: $L_1$ contains two
additional markings $z_{1,k_1+1}$ and $w_{1,k_1+1}$; there is an
$\al'$-circle $\al_{n+1}$ and a $\be'$-circle $\be_{n+1}$,
intersecting each other a point $\varrho$ with $z_{1,k_1+1}$ lying
immediately to the southeast of $\varrho$ and $w_{1,k_1+1}$ lying
immediately to the northeast of $\varrho$; and $\mc{G}$ can be obtained
from $\mc{G}'$ by deleting $z_{1,k_1+1}$, $w_{1,k_1+1}$, $\al_{n+1}$
and $\be_{n+1}$. Getting $\mc{G}'$ from $\mc{G}$ is an instance of
stabilization.

Since $\al\sbs\al'$ and $\be\sbs\be'$, $\mc{S}_{\mc{G}}$ can be
identified with the subset of $\mc{S}_{\mc{G}'}$ consisting of the
states that contain the point $\varrho$. Furthermore, since
$\CF_{\mc{G}}$ is freely generated by $\mc{S}_{\mc{G}}$ over the
$\F_2$-algebra $\mb{P}$ and $\CF_{\mc{G}'}$ is freely generated by
$\mc{S}_{\mc{G}'}$ over the $\F_2$-algebra
$\mb{P}'=\mb{P}\otimes\F_2[U_{1,k_1+1}]$, there is a natural inclusion
map $\iota\colon\CF_{\mc{G}}\hookrightarrow\CF_{\mc{G}'}$
(note, $\iota$ is just a $\mb{P}$-module map, it in general is not a
chain map). With this in mind, let us describe the naturality map
$F_{\mc{G},\mc{G}'}$.

\begin{figure}
\psfrag{p}{$\varrho$}
\begin{center}
\includegraphics[width=0.9\textwidth]{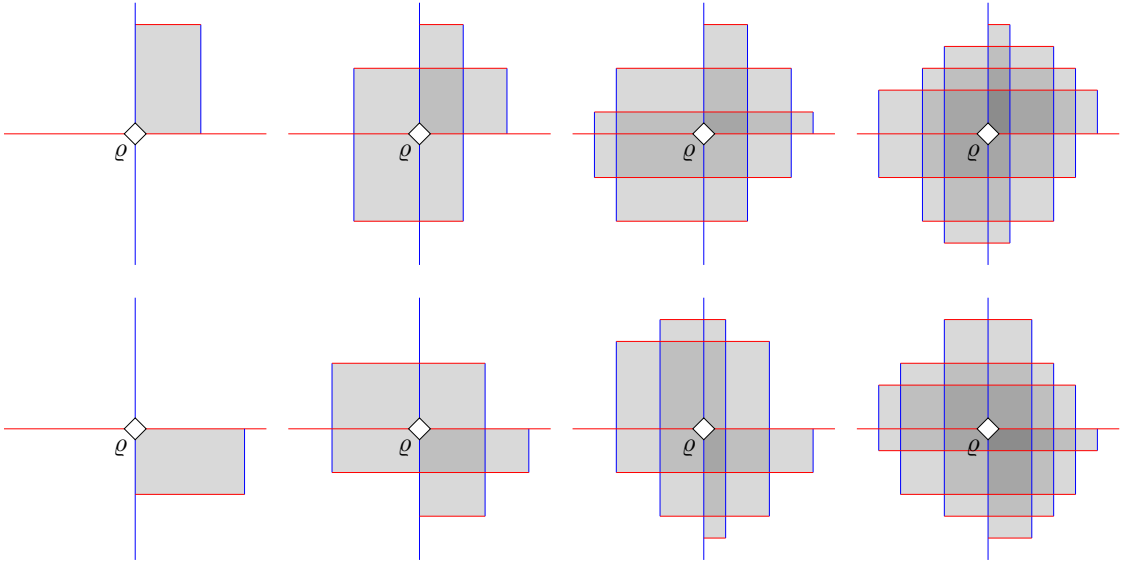}
\end{center}
\caption{Snail domains.}\label{fig:snail}
\end{figure}

For $x,y\in\mc{S}_{\mc{G}'}$, a \emph{northeast snail domain centered
at $\varrho$} is a positive domain $D\in\mc{D}_{\mc{G}'}(x,y)$ such
that the following hold: $\del D$ is an immersed circle in $T$; each
coordinate of $x$ and $y$, except possibly $\varrho$, appears with
coefficient $0$ or $\frac{1}{4}$ in $D$; finally, there is some $m\geq
0$, such that the coefficient of $D$ is $(m+1)$ immediately to the
northeast of $\varrho$, and is $m$ in the other three squares adjacent
to $\varrho$. Let $\mc{L}^1_{\varrho}(x,y)$ denote the set of all
northeast snail domain centered at $\varrho$ joining $x$ to $y$. We
have shown some elements of $\mc{L}^1_{\varrho}(x,y)$ in the first row
of \hyperref[fig:snail]{Figure \ref*{fig:snail}}.

Let us define a $\mb{P}$-module $(0,0,\ldots,0)$-graded chain map
$\wt{F}_{\mc{G},\mc{G}'}\colon\CF_{\mc{G}}\rightarrow \CF_{\mc{G}'}$
as follows. For $x\in\mc{S}_{\mc{G}}$, if $\iota(x)$ denotes the
corresponding state in $\mc{S}_{\mc{G}'}$, then
$$\wt{F}_{\mc{G},\mc{G}'}(x)=
\sum_{y\in\mc{S}_{\mc{G}'}}y\sum_{\substack{D\in\mc{L}^1_{\varrho}(\iota(x),y)\\
n_z(D)=n_{z_{1,k_1+1}}(D)}} \prod_{\substack{\imath,\jmath\\
(\imath,\jmath)\neq (1,k_1+1)}} U_{\imath,\jmath}^{n_{w_{\imath,\jmath}}(D)}.$$

There are similar maps for the other types of stabilization, which can
be defined in terms of similar looking snail domains.  All these snail
domains are very similar to the snail domains that appear
in \cite{CMPOZSzDT}; following their analysis, it is not very hard to
see that each of the stabilization maps, and in particular the above
one, is a chain map. Before we show that this map is the same as the
naturality map $F_{\mc{G},\mc{G}'}$ in
$\Mor_{K(\mc{A}_l)}(\CF_{\mc{G}},\CF_{\mc{G}'})$, let us talk
a little bit about destabilization.

Destabilization for grid diagrams is the reverse move of
stabilization. Once again, there are several variants, and once again,
we will concentrate on only one. The starting grid diagram
$\mc{G}=(T,\al,\be,z,w)$ has index $(n+1)$ and the link component
$L_1$ contains the $2k_1+2$ markings $z_{1,1},\ldots,\allowbreak
z_{1,k_1+1},w_{1,1},\allowbreak \ldots,w_{1,k_1+1}$. Furthermore,
assume that there is an $\al$-circle $\al_{n+1}$ and a $\be$-circle
$\be_{n+1}$, intersecting each other a point $\varrho$, such that
$z_{1,k_1+1}$ lies immediately to the southeast of $\varrho$ and
$w_{1,k_1+1}$ lies immediately to the southwest of $\varrho$; let
$w_{1,d}$ be the $w$-marking that lies in the same component of
$T\sm\be$ as $z_{1,k_1+1}$. The index $n$ grid diagram $\mc{G}'$ is
obtained from $\mc{G}$ by deleting $z_{1,k_1+1}$, $w_{1,k_1+1}$,
$\al_{n+1}$ and $\be_{n+1}$; this is an instance of destabilization.

Once again, there is a natural inclusion
$\iota\colon\mc{S}_{\mc{G}'}\hookrightarrow \mc{S}_{\mc{G}}$, which
induces a $\mb{P}'$-module map
$\iota\colon\CF_{\mc{G}'}\hookrightarrow \CF_{\mc{G}}$. For
$x,y\in\mc{S}_{\mc{G}}$, a \emph{southeast snail domain centered at
$\varrho$} is a positive domain $D\in\mc{D}_{\mc{G}}(x,y)$ such that
$\del D$ is an immersed circle in $T$, each coordinate of $x$ and $y$,
except possibly $\varrho$, appears with coefficient $0$ or
$\frac{1}{4}$ in $D$, and there is some $m\geq 0$, such that the
coefficient of $D$ is $(m+1)$ immediately to the southeast of
$\varrho$, and is $m$ in the other three squares adjacent to
$\varrho$. Let $\mc{L}^2_{\varrho}(x,y)$ denote the set of all
southeast snail domain centered at $\varrho$ joining $x$ to $y$. Some
elements of $\mc{L}^2_{\varrho}(x,y)$ are shown in the second row of
\hyperref[fig:snail]{Figure \ref*{fig:snail}}. \cite{CMPOZSzDT} defines a $\mb{P}'$-module
$(0,0,\ldots,0)$-graded chain map
$\wt{F}_{\mc{G},\mc{G}'}\colon\CF_{\mc{G}}\rightarrow \CF_{\mc{G}'}$
as follows: 
$$\wt{F}_{\mc{G},\mc{G}'}(U_{1,k_1+1}^m x)=
\sum_{y\in\mc{S}_{\mc{G}'}} y\sum_{\substack{D\in\mc{L}^2_{\varrho}(x,\iota(y))\\
n_z(D)=n_{z_{1,k_1+1}}(D)}} U_{1,d}^m\prod_{\imath,\jmath}
U_{\imath,\jmath}^{n_{w_{\imath,\jmath}}(D)}.$$

There are similar chain maps for the other types of
destabilization. \cite{CMPO} show that for each configuration, the
destabilization map $\wt{F}_{\mc{G},\mc{G}'}$, defined in terms of
the snail domains, is the same as the naturality map $F_{\mc{G},\mc{G}'}$
in $\Mor_{K(\mc{A}_l)}(\CF_{\mc{G}},\CF_{\mc{G}'})$. This
implies that the same is true for the stabilization maps.

\begin{lem}
If $\mc{G}'$ is obtained from $\mc{G}$ by a stabilization, then the
stabilization map $\wt{F}_{\mc{G},\mc{G}'}$, defined in terms of
the snail domains, is the same as the naturality map $F_{\mc{G},\mc{G}'}$
in $\Mor_{D(\mc{A}_l)}(\CF_{\mc{G}},\CF_{\mc{G}'})$.
\end{lem}

\begin{proof}
By naturality, we know that $F_{\mc{G},\mc{G}'}F_{\mc{G}',\mc{G}}=\Id$
in $K(\mc{A}_l)$. From \cite{CMPO} we know that
$F_{\mc{G}',\mc{G}}=\wt{F}_{\mc{G}',\mc{G}}$ in
$K(\mc{A}_l)$. Finally, it is easy to check that the chain map
$\wt{F}_{\mc{G}',\mc{G}}\wt{F}_{\mc{G},\mc{G}'}$ actually equals the
identity map on $\CF_{\mc{G}}$. Therefore,
$\wt{F}_{\mc{G},\mc{G}'}=F_{\mc{G},\mc{G}'}F_{\mc{G}',\mc{G}}\wt{F}_{\mc{G},\mc{G}'}
= F_{\mc{G},\mc{G}'}\wt{F}_{\mc{G}',\mc{G}}\wt{F}_{\mc{G},\mc{G}'}=F_{\mc{G},\mc{G}'}$. 
\end{proof}

\subsection{Renumbering}
There is yet a third type of grid move, namely, renumbering the $w$
and $z$-markings. Since we are only concerned with the link component
$L_1$, let us only consider the renumbering of the $w_1$ and the
$z_1$-markings. Let us start with a grid diagram $\mc{G}$ with $k_1$
$w_1$-markings and $k_1$ $z_1$-markings. Fix $\si,\tau\in\mf{S}_{k_1}$. Let
$\mc{G}'$ be the grid diagram obtained from $\mc{G}$ by renaming
$w_{1,i}$ as $w_{1,\si(i)}$ and $z_{1,i}$ as $z_{1,\tau(i)}$ for $1\leq
i\leq k_1$; furthermore, if $w_{1,s_1}$ is the special $w_1$-marking
in $\mc{G}$, then $w_{1,\si(s_1)}$ is the special $w'_1$-marking in
$\mc{G}'$.

The grid diagrams $\mc{G}$ and $\mc{G}'$ represent the same pointed
link. The naturality map $F_{\mc{G},\mc{G}'}$ is $U_{i,j}$-equivariant
for $i\neq 1$, and sends $x\prod_j U_{1,j}^{m_j}$ to $x\prod_j
U_{1,\si(j)}^{m_j}$. This is simply because the $U_{i,j}$ variables
are not indexed by pairs of integers $i,j$, but are rather indexed by
the $w$-markings themselves. Indeed, we should have written $U_{i,j}$
as $U_{w_{i,j}}$. Precise notations lead to triple subscripts
(as in $U_{w_{i,s_i}}$!),  so we have chosen to avoid them.

\subsection{Changing the special marking}\label{subsec:changesp}
 Let $\mc{G}=\mc{G}_1$ be a grid
diagram for the pointed link $(L,p)$ such that $s_1=1$ (i.e.\  the
basepoint $p_1\in L_1$ is represented by the $w$-marking $w_{1,1}$),
and the markings that appear in $L_1$ are, in order,
$w_{1,1},z_{1,1},\ldots,w_{1,k_1},z_{1,k_1}$ with $k_1>1$. For $2\leq j\leq k_1$,
let $\mc{G}_j$ be the grid diagram where $s_1=j$, but which is
otherwise identical to $\mc{G}$. Observe that the grid diagrams
$\mc{G}_j$ represent the same link, but not the same pointed link. The
chain complexes $\CF_{\mc{G}_j}$ are identical as
$\F_2[U_2,\ldots,U_l]$-modules; the $U_1$-action is multiplication by
$U_{1,j}$.

For distinct $j,j'\in\{1,\ldots,k_1\}$, we will define
$U_i$-equivariant chain maps $f_{j,j'}$ from $\CF_{\mc{G}_j}$ to
$\CF_{\mc{G}_{j'}}$. Define
$c_{j,j'}\colon\CF_{\mc{G}_j}\rightarrow \CF_{\mc{G}_{j'}}$ as follows: it
is $U_{\imath,\jmath}$-equivariant for $(\imath,\jmath)\neq (1,j)$ and
sends $U_{1,j}^m x$ to $U_{1,j'}^m x$. It is clearly an
$\F_2[U_1,\ldots,U_l]$-module map, but in general is not a chain
map. Let $J(j,j')\sbseq\{1,\ldots,k_1\}$ be defined as: $\jmath\in
J(j,j')$ if and only if $z_{1,\jmath}$ appears in the arc joining
$w_{1,j}$ to $w_{1,j'}$ in the oriented link component $L_1$. Define
$$f_{j,j'}=c_{j,j'}\Phi_{1,j}\sum_{\jmath\in
J(j,j')}\Psi_{1,\jmath}.$$

\begin{thm}\label{thm:changespecial}
The map $f_{j,j'}$, as defined above, is an
$\F_2[U_1,\ldots,U_l]$-module chain map from $\CF_{\mc{G}_j}$ to
$\CF_{\mc{G}_{j'}}$. Furthermore,
$f_{k_1,1}f_{1,k_1}=\Id+\Phi_1\Psi_1=\Id+\Psi_1\Phi_1$, and for $2\leq
j\leq k_1-1$, $f_{j,j+1}f_{1,j}=f_{1,j+1}$ in $K(\mc{A}_l)$
\end{thm}

\begin{proof}
The maps $\Phi_{i,j}$ and $\Psi_{i,j}$ are
$U_{\imath,\jmath}$-equivariant. The map $c_{j,j'}$ is also
$U_i$-equivariant since $[c_{j,j'}\colon U_1]=c_{j,j'}U_1+U_1
c_{j,j'}=c_{j,j'}U_{1,j}+U_{1,j'}c_{j,j'}=0$. Therefore, the maps
$f_{j,j'}$ are $U_i$-equivariant.

To see that $f_{j,j'}$ is a chain map, recall the commutator
relations: $[\Phi_{1,\jmath}\colon\del]=0$ and
$[\Psi_{1,\jmath}\colon\del]=U_{1,\jmath}+U_{1,\jmath+1}$ (the second
index being numbered modulo $k_1$). It is also easy to see that
$[c_{j,j'}\colon\del]=(U_{1,j}+U_{1,j'})c_{j,j'}\Phi_{1,j}$ (since
$w_{1,j}$ can appear in a rectangle at most once). Therefore,
\begin{align*}
[f_{j,j'}\colon\del]&=([c_{j,j'}\colon\del]\Phi_{1,j}\sum_{\jmath\in
J(j,j')}\Psi_{1,\jmath})+(c_{j,j'}[\Phi_{1,j}\colon\del]\sum_{\jmath\in
J(j,j')}\Psi_{1,\jmath})+(c_{j,j'}\Phi_{1,j}\sum_{\jmath\in
J(j,j')}[\Psi_{1,\jmath}\colon\del])\\
&=((U_{1,j}+U_{1,j'})c_{j,j'}\Phi^2_{1,j}\sum_{\jmath\in
J(j,j')}\Psi_{1,\jmath})+(c_{j,j'}\Phi_{1,j}(U_{1,j}+U_{1,j'}))\\
&=0+c_{j,j'}(U_{1,j}+U_{1,j'})\Phi_{1,j}\qquad\text{(since $\Phi^2_{1,j}=0$
in grid diagrams)}\\
&=0.
\end{align*}

For the second part of the theorem, let $\Id_{j,j'}$ be the identity
map from $\CF_{\mc{G}_j}$ to $\CF_{\mc{G}_{j'}}$. It is not
$U_1$-equivariant; indeed, $[\Id_{j,j'}\colon
U_1]=(U_{1,j}+U_{1,j'})\Id_{j,j'}$. Define
$K_{j,j'}\colon \CF_{\mc{G}_j}\rightarrow\CF_{\mc{G}_{j'}}$ as
follows: it is $U_{\imath,\jmath}$-equivariant for
$(\imath,\jmath)\neq (1,j)$ and sends $U_{1,j}^m x$ to
$\frac{U^m_j+U^m_{j'}}{U_j+U_{j'}}x$. It is easy to see that
$[K_{j,j'}\colon U_1] = U_{1,j'}K_{j,j'}+K_{j,j'}U_{1,j} = \Id_{j,j'}$;
furthermore, $[K_{j,j'}\colon\del]= c_{j,j'}\Phi_{1,j}$.

We have
\begin{align*}
f_{k_1,1}f_{1,k_1}+\Phi_1\Psi_1&=c_{k_1,1}\Phi_{1,k_1}\Psi_{1,k_1}
c_{1,k_1}\Phi_{1,1}(\Psi_{1,1}+\cdots+\Psi_{1,k_1-1})+\Phi_{1,1}((\Psi_{1,1}+\cdots+\Psi_{1,k_1})\\
&=(c_{k_1,1}\Phi_{1,k_1}\Psi_{1,k_1}
\Id_{1,k_1}+\Id)c_{1,1}\Phi_{1,1}(\Psi_{1,1}+\cdots+\Psi_{1,k_1-1})+c_{1,1}\Phi_{1,1}\Psi_{1,k_1}
\end{align*}
and for $2\leq j\leq k_1-1$,
\begin{align*}
f_{1,j}f_{j,j+1}+f_{1,j+1}&=c_{j,j+1}\Phi_{1,j}\Psi_{1,j}
c_{1,j}\Phi_{1,1}(\Psi_{1,1}+\cdots+\Psi_{1,j-1})+c_{1,j+1}\Phi_{1,1}((\Psi_{1,1}+\cdots+\Psi_{1,j})\\
&=(c_{j,j+1}\Phi_{1,j}\Psi_{1,j}
\Id_{j+1,j}+\Id)c_{1,j+1}\Phi_{1,1}(\Psi_{1,1}+\cdots+\Psi_{1,j-1})+c_{1,j+1}\Phi_{1,1}\Psi_{1,j}.
\end{align*}
Therefore, for $2\leq j\leq k_1$, we are interested in the map $g_j=(c_{j,j+1}\Phi_{1,j}\Psi_{1,j}
\Id_{j+1,j}+\Id)c_{1,j+1}\Phi_{1,1}(\Psi_{1,1}+\cdots+\Psi_{1,j-1})+c_{1,j+1}\Phi_{1,1}\Psi_{1,j}$. We
want to show that $g_{k_1}=\Id$ and $g_j=0$ for $j\leq k_1-1$ in $K(\mc{A}_l)$.

Consider the map
$H_j=K_{j,j+1}\Psi_{1,j}\Id_{j+1,j}c_{1,j+1}\Phi_{1,1}(\Psi_{1,1}+\cdots+\Psi_{1,j-1})$. We
have
\begin{align*}
[H_j\colon U_1]&=[K_{j,j+1}\colon
U_1]\Psi_{1,j}\Id_{j+1,j}c_{1,j+1}\Phi_{1,1}(\Psi_{1,1}+\cdots+\Psi_{1,j-1})\\
 &\qquad {}+K_{j,j+1}\Psi_{1,j}
[\Id_{j+1,j}\colon
U_1]c_{1,j+1}\Phi_{1,1}(\Psi_{1,1}+\cdots+\Psi_{1,j-1})\\
&=\Id_{j,j+1}\Psi_{1,j}\Id_{j+1,j}c_{1,j+1}\Phi_{1,1}(\Psi_{1,1}+\cdots+\Psi_{1,j-1})\\
&\qquad {}+K_{j,j+1}\Psi_{1,j}
(U_{1,j}+U_{1,j+1})\Id_{j+1,j}c_{1,j+1}\Phi_{1,1}(\Psi_{1,1}+\cdots+\Psi_{1,j-1})\\
&=\Id_{j,j+1}\Psi_{1,j}\Id_{j+1,j}c_{1,j+1}\Phi_{1,1}(\Psi_{1,1}+\cdots+\Psi_{1,j-1})\\
&\qquad
{}+\Id_{j,j+1}\Psi_{1,j}\Id_{j+1,j}c_{1,j+1}\Phi_{1,1}(\Psi_{1,1}+\cdots+\Psi_{1,j-1})\\
&=0.
\end{align*}
Therefore, $H_j$ is a $U_i$-equivariant map from $\CF_{\mc{G}}$ to
$\CF_{\mc{G}_{j+1}}$. We also have 
\begin{align*}
[H_j\colon \del]&=[K_{j,j+1}\colon
\del]\Psi_{1,j}\Id_{j+1,j}c_{1,j+1}\Phi_{1,1}(\Psi_{1,1}+\cdots+\Psi_{1,j-1})\\
&\qquad {}+
K_{j,j+1}[\Psi_{1,j}\colon\del]\Id_{j+1,j}c_{1,j+1}\Phi_{1,1}(\Psi_{1,1}+\cdots+\Psi_{1,j-1})\\
&\qquad {}+
K_{j,j+1}\Psi_{1,j}\Id_{j+1,j}[c_{1,j+1}\colon\del]\Phi_{1,1}(\Psi_{1,1}+\cdots+\Psi_{1,j-1})\\
&\qquad {}+K_{j,j+1}\Psi_{1,j}\Id_{j+1,j}c_{1,j+1}\Phi_{1,1}([\Psi_{1,1}\colon\del]+\cdots+
[\Psi_{1,j-1}\colon\del])\\
&=c_{j,j+1}\Phi_{1,j}\Psi_{1,j}\Id_{j+1,j}c_{1,j+1}\Phi_{1,1}(\Psi_{1,1}+\cdots+\Psi_{1,j-1})\\
&\qquad {}+
K_{j,j+1}(U_{1,j}+U_{1,j+1})\Id_{j+1,j}c_{1,j+1}\Phi_{1,1}(\Psi_{1,1}+\cdots+\Psi_{1,j-1})\\
&\qquad {}+
K_{j,j+1}\Psi_{1,j}\Id_{j+1,j}(U_{1,1}+U_{1,j+1})c_{1,j+1}\Phi^2_{1,1}(\Psi_{1,1}+\cdots+\Psi_{1,j-1})\\
&\qquad
{}+K_{j,j+1}\Psi_{1,j}\Id_{j+1,j}c_{1,j+1}\Phi_{1,1}(U_{1,1}+U_{1,j})\\
&=c_{j,j+1}\Phi_{1,j}\Psi_{1,j}\Id_{j+1,j}c_{1,j+1}\Phi_{1,1}(\Psi_{1,1}+\cdots+\Psi_{1,j-1})\\
&\qquad {}+
\Id_{j,j+1}\Id_{j+1,j}c_{1,j+1}\Phi_{1,1}(\Psi_{1,1}+\cdots+\Psi_{1,j-1})\\
&\qquad
{}+0+K_{j,j+1}(U_{1,j+1}+U_{1,j})\Psi_{1,j}\Id_{j+1,j}c_{1,j+1}\Phi_{1,1}\\
&=(c_{j,j+1}\Phi_{1,j}\Psi_{1,j}
\Id_{j+1,j}+\Id)c_{1,j+1}\Phi_{1,1}(\Psi_{1,1}+\cdots+\Psi_{1,j-1})+\Psi_{1,j}c_{1,j+1}\Phi_{1,1}.
\end{align*}
Thus,
$[H_j\colon\del]+g_j=\Psi_{1,j}c_{1,j+1}\Phi_{1,1}+c_{1,j+1}\Phi_{1,1}\Psi_{1,j}$. We
want to show that it is $\Id$ for $j=k_1$, and is zero for $j\leq
k_1-1$ in $K(\mc{A}_l)$.

For $j\leq k_1-1$, consider the map
$K_{1,j+1}\Psi_{1,j}+\Psi_{1,j}K_{1,j+1}$. It is straightforward to
check that it is $U_i$-equivariant and its commutator with $\del$ is
$\Psi_{1,j}c_{1,j+1}\Phi_{1,1}+c_{1,j+1}\Phi_{1,1}\Psi_{1,j}$.

For the last remaining case, consider the $U_i$-equivariant map
$H\colon \CF_{\mc{G}}\rightarrow\CF_{\mc{G}}$, defined as follows:
$$H(x)=\sum_{y\in\mc{S}_{\mc{G}}}y\sum_{\substack{R\in\mc{R}_{\mc{G}}(x,y)\\n_{w_{1,1}}(R)=1\\
n_z(R)=n_{z_{1,k_1}}(R)=1}}\prod_{(\imath,\jmath)\neq
(1,1)}U_{\imath,\jmath}^{n_{w_{\imath,\jmath}}(R)}.$$
It is easy to see that, since $k_1>1$,
$[H\colon\del]=\Phi_{1,1}\Psi_{1,k_1}+\Psi_{1,k_1}\Phi_{1,1}+\Id$. This
concludes the proof.
\end{proof}

\subsection{The main theorem}
At this point, modulo the technical lemma
\hyperref[lem:technical]{Lemma \ref*{lem:technical}}, we are ready to
prove \hyperref[thm:main]{Theorem
  \ref*{thm:main}}. \hyperref[lem:technical]{Lemma
  \ref*{lem:technical}} is hard to motivate on its own, so let us
postpone it until we need it.

\begin{proof}[Proof of Theorem \ref*{thm:main}]
As mentioned earlier, we are proving this for $i=1$. Due
to \hyperref[thm:invariance]{Theorem \ref*{thm:invariance}}, we are
free to choose our Heegaard diagram. So let $\mc{H}$ be an index $n$
grid diagram representing the pointed link $(L,p)$; the Heegaard
surface $\Si$ is a torus, so let us call it $T$; let $J_s$ be a
constant path of almost complex structure on $\Sym^n(T)$ induced from
some complex structure $\mf{j}$ on $T$; as usual, let
$w_{1,1},z_{1,1},\ldots,w_{1,k_1},z_{1,k_1}$ be the markings, in
order, on the oriented link component $L_1$, with $k_1>1$ and
$w_{1,1}$ being the special $w_1$-marking. Later on, during the proof
of \hyperref[lem:technical]{Lemma
\ref*{lem:technical}}, we will impose further restrictions on
$\mc{H}$; but for now, let us not be concerned with such restrictions.

As in \hyperref[subsec:changesp]{Subsection
\ref*{subsec:changesp}}, for $1\leq j\leq k_1$, let $\mc{G}_j$ be the
grid diagram where $w_{1,j}$ is the special $w_1$-marking, but which
is otherwise identical to $\mc{H}$; let $(L,p^{(j)})$ be the pointed
link that $\mc{G}_j$ represents (in particular, $\mc{G}_1=\mc{H}$ and
$p^{(1)}=p$). For each $j$, starting at $(\mc{G}_j,J_s)$, we will do a
stabilization, a sequence of horizontal commutations, a sequence of
vertical commutations, a renumbering, a destabilization and finally a
change of complex structure, to obtain a diagram
$(\mc{G}'_{j+1},J'_s)$; we will then apply a self-homeomorphism
$\tau_{j,j+1}$ of $T$ that sends $\mc{G}'_{j+1}$ to $\mc{G}_{j+1}$ and
$J'_s$ to $J_s$ (the relevant subscripts are numbered modulo
$k_1$). This process is best described by \hyperref[fig:themove]{Figure
\ref*{fig:themove}} (the shaded region might contain additional $w$
and $z$-markings).

\begin{figure}
\captionsetup[subfloat]{font=normalsize,labelformat=empty,width={0.45\textwidth}}
\begin{center}
\subfloat[We start with the
diagram $(\mc{G}_{j},J_s).$]{\psfrag{w1}{$w_{1,j+1}$}\psfrag{w2}{$w_{1,j}$}\psfrag{w3}{}\psfrag{z1}{$z_{1,j}$}\psfrag{z2}{}\includegraphics[width=0.3\textwidth]{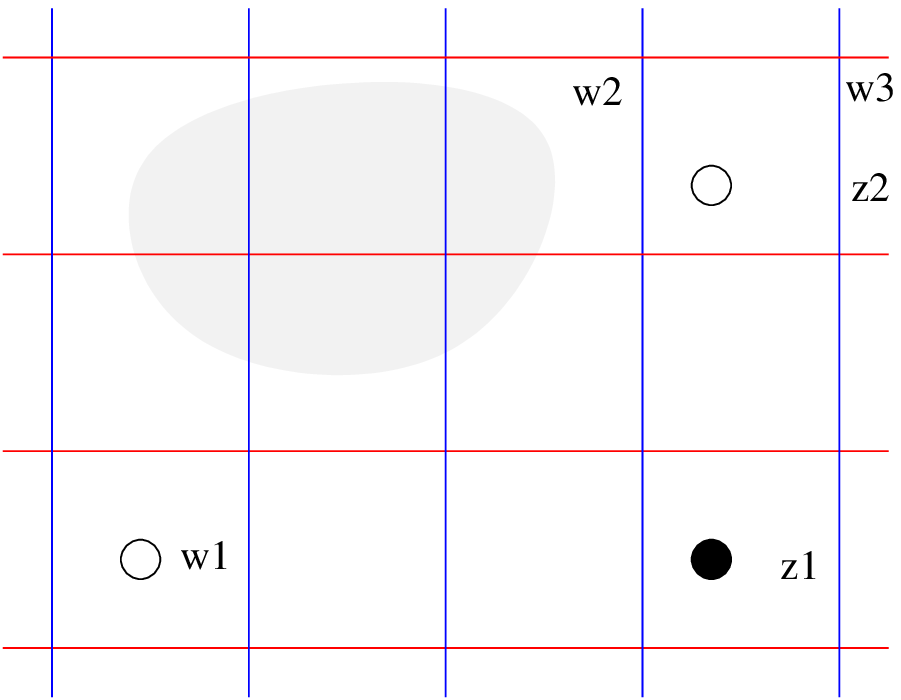}}
\hspace{0.2\textwidth}\subfloat[We stabilize by adding $w_{1,k_1+1}$ and $z_{1,k_1+1}$.]{\psfrag{w1}{$w_{1,j+1}$}\psfrag{w2}{$w_{1,j}$}\psfrag{w3}{$w_{1,k_1+1}$}\psfrag{z1}{$z_{1,j}$}\psfrag{z2}{$z_{1,k_1+1}$}\includegraphics[width=0.3\textwidth]{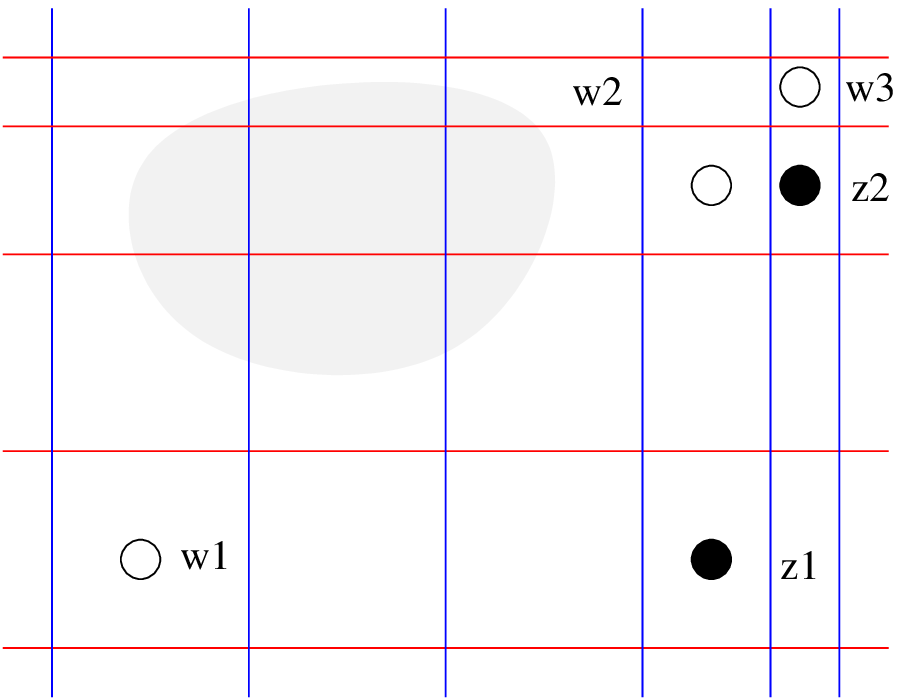}}\\
\subfloat[We do a sequence of horizontal commutations.]{\psfrag{w1}{$w_{1,j+1}$}\psfrag{w2}{$w_{1,j}$}\psfrag{w3}{$w_{1,k_1+1}$}\psfrag{z1}{$z_{1,j}$}\psfrag{z2}{$z_{1,k_1+1}$}\includegraphics[width=0.3\textwidth]{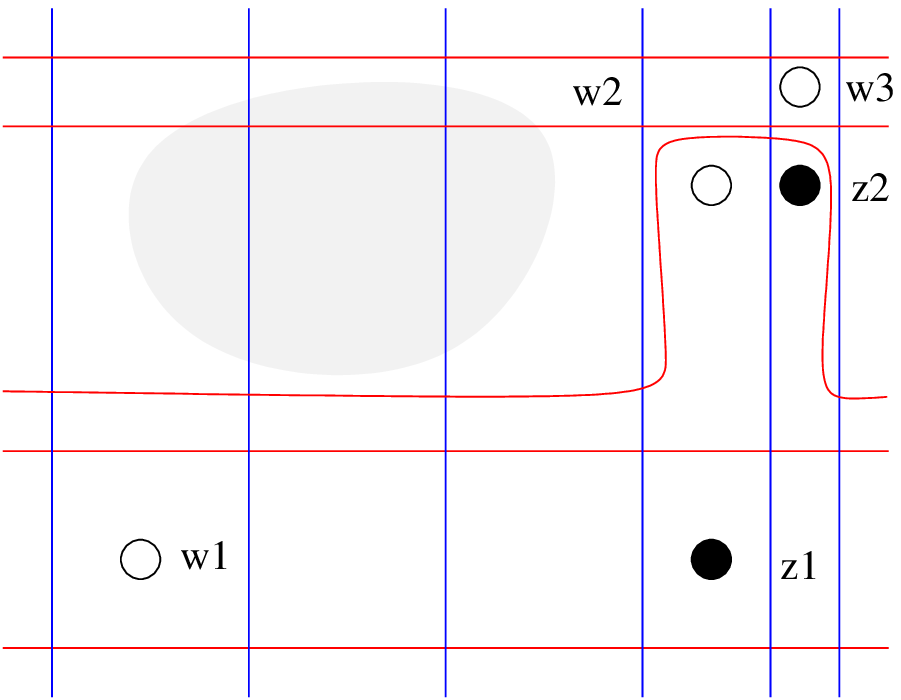}}\hspace{0.2\textwidth}
\subfloat[We do a sequence of vertical commutations.]{\psfrag{w1}{$w_{1,j+1}$}\psfrag{w2}{$w_{1,j}$}\psfrag{w3}{$w_{1,k_1+1}$}\psfrag{z1}{$z_{1,j}$}\psfrag{z2}{$z_{1,k_1+1}$}\includegraphics[width=0.3\textwidth]{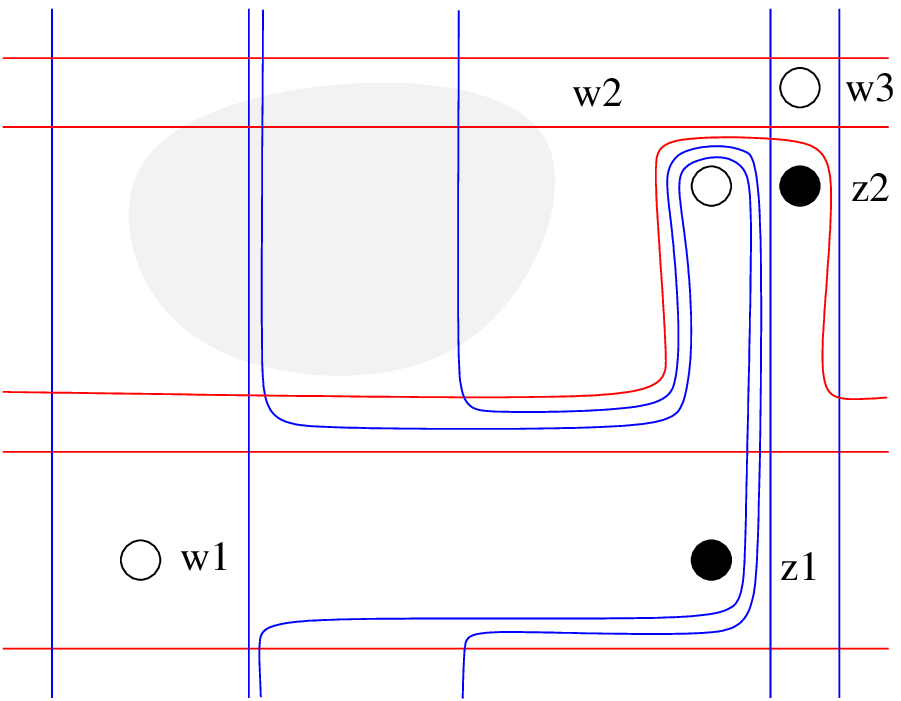}}\\
\subfloat[We renumber the $w$ and $z$-markings.]{\psfrag{w1}{$w_{1,k_1+1}$}\psfrag{w2}{$w_{1,j+1}$}\psfrag{w3}{$w_{1,j}$}\psfrag{z1}{$z_{1,k_1+1}$}\psfrag{z2}{$z_{1,j}$}\includegraphics[width=0.3\textwidth]{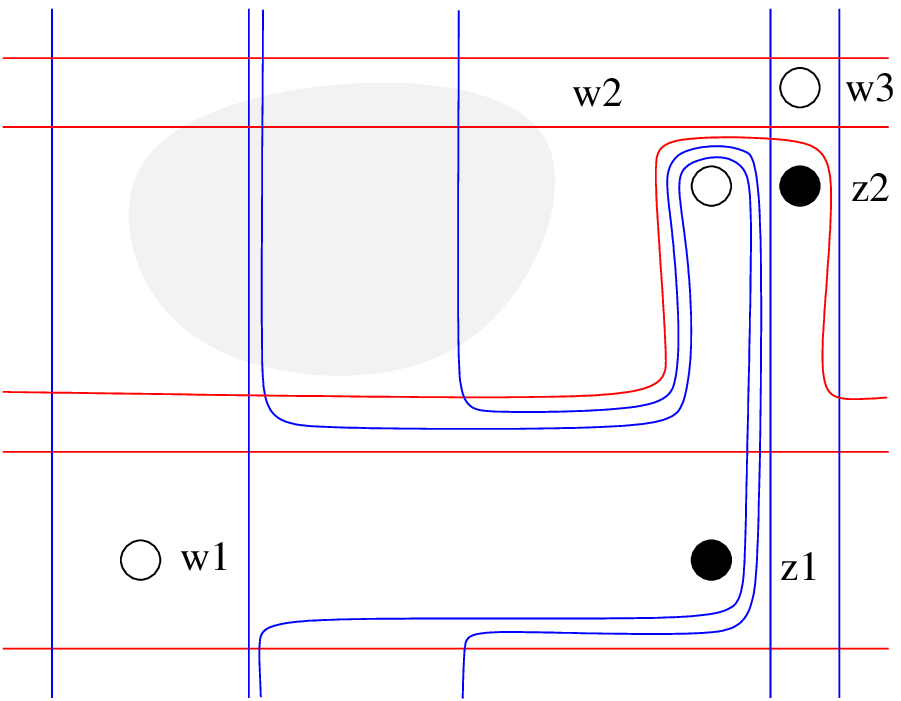}}\hspace{0.2\textwidth}
\subfloat[We destabilize by deleting $w_{1,k_1+1}$ and
$z_{1,k_1+1}$. There is an obvious homeomorphism $\tau_{j,j+1}$ which
converts this diagram to
$\mc{G}_{j+1}$.]{\psfrag{w1}{}\psfrag{w2}{$w_{1,j+1}$}\psfrag{w3}{$w_{1,j}$}\psfrag{z1}{}\psfrag{z2}{$z_{1,j}$}\includegraphics[width=0.3\textwidth]{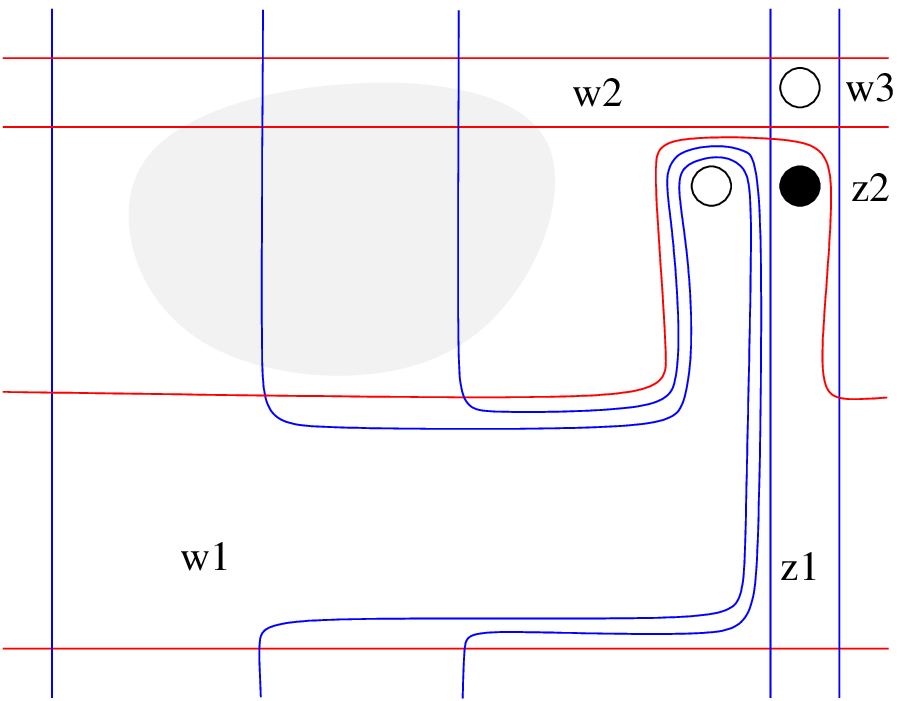}}\\
\subfloat[We change the complex structure on $T$ from $\mf{j}$ to
$\mf{j}'=\tau_{j,j+1}^{-1}(\mf{j})$ to get the diagram
$(\mc{G}'_{j+1},J'_s)$.]{\psfrag{w1}{}\psfrag{w2}{$w_{1,j+1}$}\psfrag{w3}{$w_{1,j}$}\psfrag{z1}{}\psfrag{z2}{$z_{1,j}$}\includegraphics[width=0.3\textwidth]{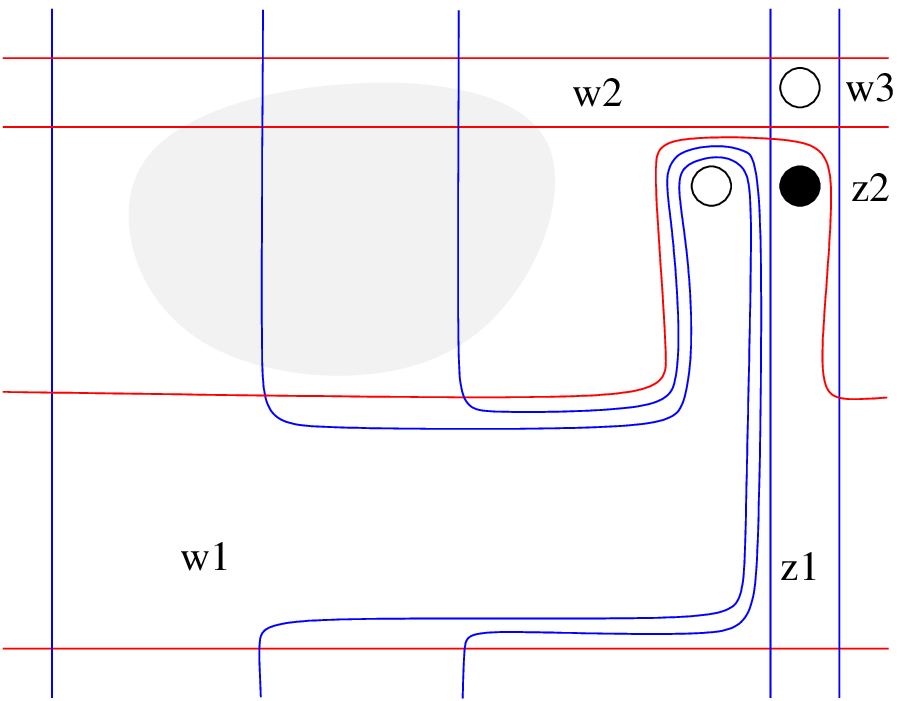}}\hspace{0.2\textwidth}
\subfloat[We apply the homeomorphism $\tau_{j,j+1}$ to get the final
diagram $(\mc{G}_{j+1},J_s)$.]{\psfrag{w1}{$w_{1,j+1}$}\psfrag{w2}{$w_{1,j}$}\psfrag{w3}{}\psfrag{z1}{$z_{1,j}$}\psfrag{z2}{}\includegraphics[width=0.3\textwidth]{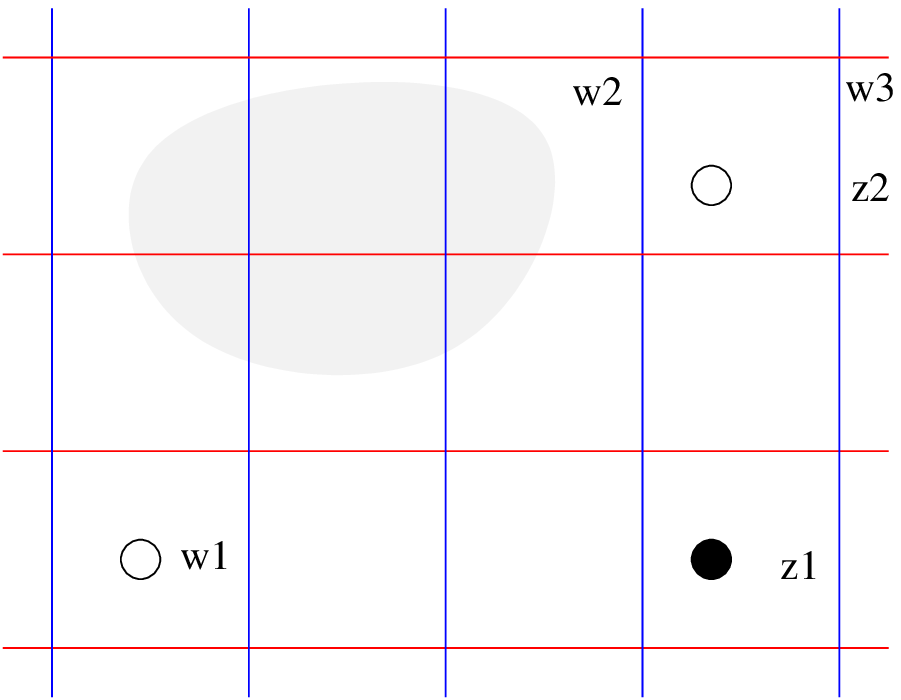}}
\end{center}
\caption{The sequence of moves converting $(\mc{G}_j,J_s)$ to $(\mc{G}_{j+1},J_s)$.}\label{fig:themove}
\end{figure}


Let us now discuss a minor subtlety that we had held back for this
long. If we stabilize or destabilize a grid diagram (or a Heegaard
diagram in general), the new diagram does not represent the same
pointed link; it represents an isotopic link. However, we want the two
diagrams to represent `nearby' links. Towards this end, let us broaden the class of links that
a Heegaard diagram represents.

Given a grid diagram $\mc{G}=(T,\al,\be,z,w)$, let $U_{\al}$ and
$U_{\be}$ be the handlebodies specified by the data $(T,\al)$ and
$(T,\be)$, respectively. To obtain \emph{a} link that $\mc{G}$
represents, in each component of $\Si\sm\al$ join the $z$-marking to
the $w$-marking by an embedded path, and push the interior of the path
in the the interior of $U_{\al}$, and in each component of $\Si\sm\be$
join the $w$-marking to the $z$-marking by an embedded path, and push
the interior of the path in the the interior of $U_{\be}$. If these
`pushes' are small, then each of the links that $\mc{G}$ represents
is supported in a small neighborhood of $T$. From now on, let us
assume that the links represented by grid diagrams are supported in a
small neighborhood of the Heegaard torus.

Therefore, when we stabilize a grid diagram $\mc{G}$ to get a grid
diagram $\mc{G}'$, even though the process changes the underlying
links, the change is not drastic; the process merely introduces a
`kink', cf.\ \hyperref[fig:kink]{Figure \ref*{fig:kink}}. Therefore,
modulo introducing a small kink, and then removing another, the
self-homeomorphism $\tau_{j,j+1}$ of $T$ induces a self-homeomorphism
$\wt{\tau}_{j,j+1}$ of $(S^3,L)$; it is constant outside a
small neighborhood of the oriented arc that joins $p^{(j)}_1=w_{1,j}$ to
$p^{(j+1)}_1=w_{1,j+1}$ on $L_1$, and it sends $p^{(j)}_1$ to
$p^{(j+1)}_1$. Thus, the composition
$\wt{\tau}_{k_1,1}\cdots\wt{\tau}_{2,3}\wt{\tau}_{1,2}$ induces the
mapping class group element $\si_1\in\MCG(S^3,L,p)$.

\begin{figure}
\begin{center}
\includegraphics[width=0.6\textwidth]{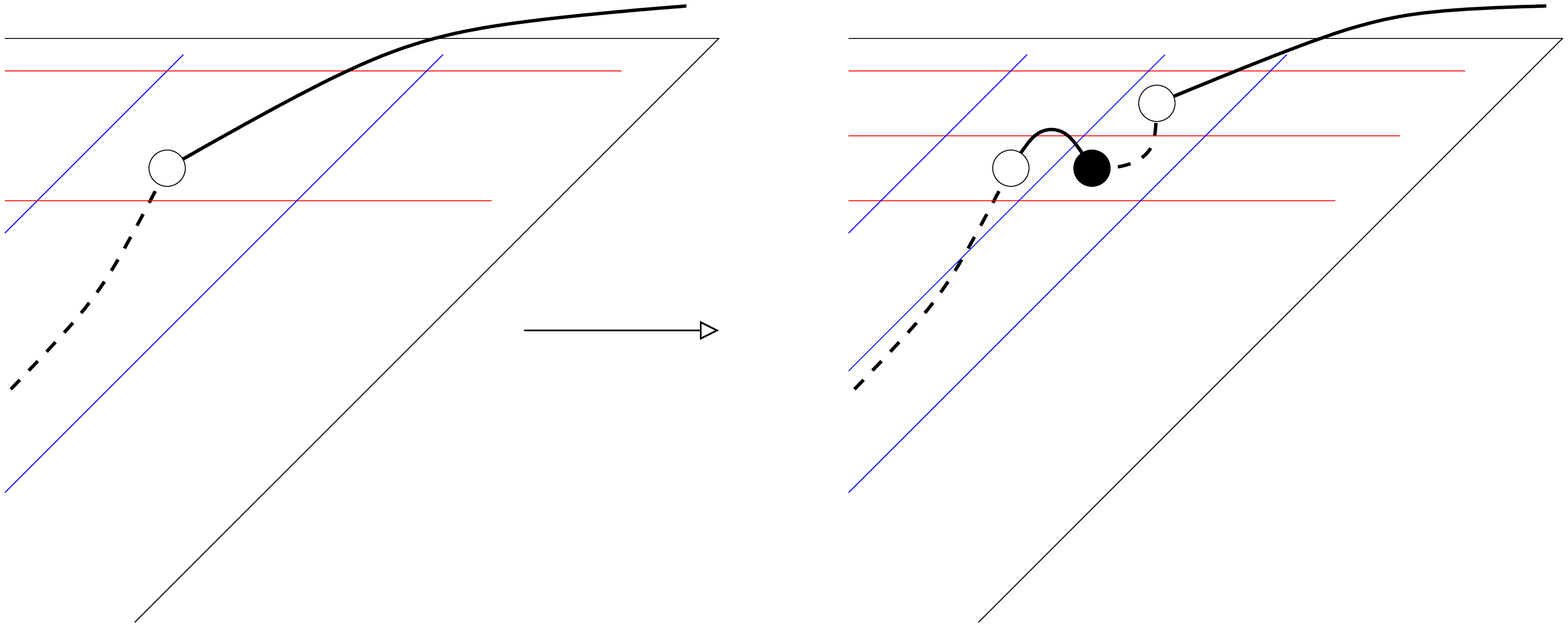}
\end{center}
\caption{Stabilization introduces a kink in the link.}\label{fig:kink}
\end{figure}

The sequence of moves that converts $(\mc{G}_j, J_s)$ to
$(\mc{G}'_{j+1},J'_s)$ induces a map from $\CF_{(\mc{G}_j,J_s)}$ to
$\CF_{(\mc{G}'_{j+1},J'_s)}$ (the map is the composition of a
stabilization map, some horizontal commutation maps, some vertical
commutation maps, a renumbering map, a destabilization map and a
change of complex structure map).  The automorphism $\tau_{j,j+1}$
induces an identification between $\CF_{(\mc{G}'_{j+1},J'_s)}$ and
$\CF_{(\mc{G}_{j+1},J_s)}$. The composition is an
$\F_2[U_1,\ldots,U_l]$-module chain map from $\CF_{\mc{G}_j}$ to
$\CF_{\mc{G}_{j+1}}$. \hyperref[lem:technical]{Lemma
  \ref*{lem:technical}} will prove that, for certain types of grid
diagrams called superstabilized grid diagrams, this map is in
fact the map $f_{j,j+1}=c_{j,j+1}\Phi_{1,j}\sum_{\jmath\in
  J(j,j+1)}\Psi_{1,\jmath}=c_{j,j+1}\Phi_{1,j}\Psi_{1,j}$.

A \emph{superstabilized grid diagram} is a grid diagram where every
$w$ or $z$-marking lies in a $2\times 2$ square $S$ which looks like
one of the following: $S$ contains a $z$-marking at the southeast
square and two $w$-markings, one each at the northeast square and the
southwest square; or $S$ contains a $w$-marking at the northwest
square and two $z$-markings, one each at the northeast square and the
southwest square. It is clear that any grid diagram can be converted
to a superstabilized one by stabilizing sufficiently many
times. Therefore, we might assume that our original grid diagram
$\mc{H}$ was superstabilized to start with.

\begin{lem}\label{lem:technical}
  Assuming that $\mc{G}_j$ is a superstabilized grid diagram, the
  sequence of moves from \hyperref[fig:themove]{Figure
    \ref*{fig:themove}}, which converts $(\mc{G}_j, J_s)$ to
  $(\mc{G}_{j+1},J_s)$, induces the map
  $f_{j,j+1}=c_{j,j+1}\Phi_{1,j}\Psi_{1,j}$ from $\CF_{\mc{G}_j}$ to
  $CF_{\mc{G}'_j}$.
\end{lem}

\begin{proof}
  Let $\mc{G}^m$ be the grid diagram shown in the $m\ith$ picture of
  \hyperref[fig:themove]{Figure \ref*{fig:themove}} (hence,
  $\mc{G}^0=\mc{G}_j$ and $\mc{G}^8=\mc{G}_{j+1}$). The map induced
  by the sequence of moves is the composition of the following maps: a
  stabilization map $s$ from $\CF_{\mc{G}^1}$ to $\CF_{\mc{G}^2}$; the
  composition $h$ of some horizontal commutation maps, which maps from
  $\CF_{\mc{G}^2}$ to $\CF_{\mc{G}^3}$; the composition $v$ of some
  vertical commutation maps, which maps from $\CF_{\mc{G}^3}$ to
  $\CF_{\mc{G}^4}$; a renumbering map $r$ from $\CF_{\mc{G}^4}$ to
  $\CF_{\mc{G}^5}$; a destabilization map $d$ from $\CF_{\mc{G}^5}$ to
  $\CF_{\mc{G}^6}$; a change of complex structure map from
  $\CF_{\mc{G}^6}$ to $\CF_{\mc{G}^7}=\CF_{\mc{G}^6}$, which is the
  identity map; and finally an identification map induced by
  $\tau_{j,j+1}$ from $\CF_{\mc{G}^7}$ to $\CF_{\mc{G}^8}$. Therefore,
  if we identify $\CF_{\mc{G}^6}$ with $\CF_{\mc{G}^1}$ by an
  identification map $\iota$, the map in
  question is simply the composition $\iota drvhs$.

\begin{figure}
\psfrag{wj}{$w_{1,j}$}
\psfrag{zj}{$z_{1,j}$}
\psfrag{wj1}{$w_{1,j+1}$}
\psfrag{zj1}{$z_{1,j+1}$}
\psfrag{wja}{$w_{1,j-1}$}
\psfrag{zja}{$z_{1,j-1}$}
\psfrag{wj2}{$w_{1,j+2}$}
\begin{center}
\includegraphics[width=0.8\textwidth]{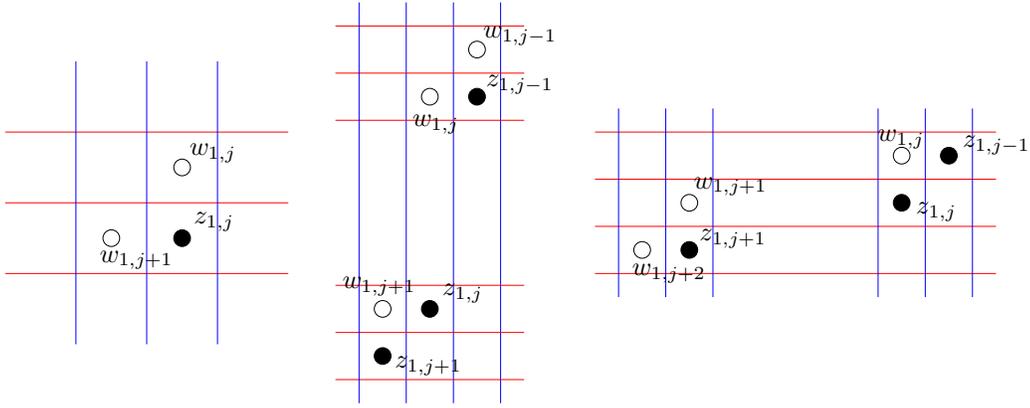}
\end{center}
\caption{The three possible configurations of $\mc{G}_j$.}\label{fig:3config}
\end{figure}

  Since $\mc{G}_j$ is a superstabilized grid diagram, it actually
  looks like one of the three configurations of
  \hyperref[fig:3config]{Figure \ref*{fig:3config}}. The first case is
  the easiest, so we will do it for warm-up. The second and the third
  case are similar, so we will only do the second case.

\begin{figure}
\psfrag{wj}{$w_{1,j}$}
\psfrag{zj}{$z_{1,j}$}
\psfrag{wj1}{$w_{1,j+1}$}
\psfrag{zk}{$z_{1,k_1+1}$}
\psfrag{wk}{$w_{1,k_1+1}$}
\psfrag{s}{$s$}
\psfrag{r}{$r$}
\psfrag{d}{$d$}
\psfrag{g1}{$\mc{G}^1$}
\psfrag{g2}{$\mc{G}^2=\mc{G}^3=\mc{G}^4$}
\psfrag{g3}{$\mc{G}^5$}
\psfrag{g4}{$\mc{G}^6$}
\begin{center}
\includegraphics[width=0.9\textwidth]{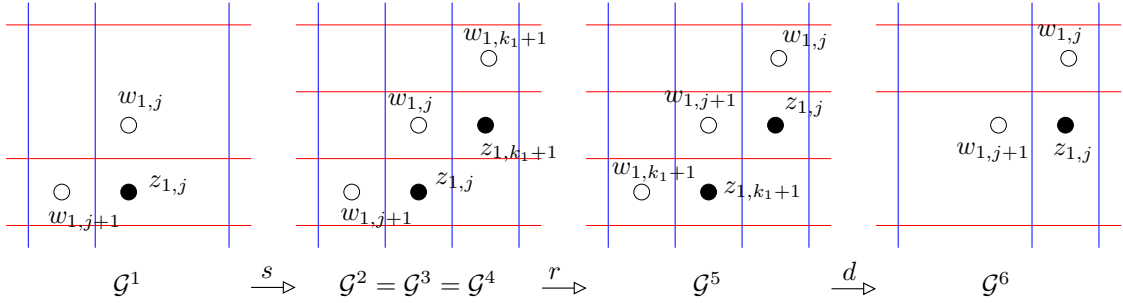}
\end{center}
\caption{The sequence of moves in the first case.}\label{fig:case1}
\end{figure}

  In the first case, $\mc{G}^2=\mc{G}^3=\mc{G}^4$, and $h=v=\Id$. The
  sequence of moves is shown in \hyperref[fig:case1]{Figure
    \ref*{fig:case1}}. We want to show that the composition $\iota
  drs$ equals the map $c_{j,j+1}\Phi_{1,j}\Psi_{1,j}$. It is easy to
  see that the part $c_{j,j+1}$ comes from renumbering map $r$ and the
  destabilization map $d$. We only have to check that the domains that
  contribute to the map $drs$ correspond to the domains that contribute
  to the map $\Phi_{1,j}\Psi_{1,j}$.

\begin{figure}
\begin{center}
\includegraphics[width=\textwidth]{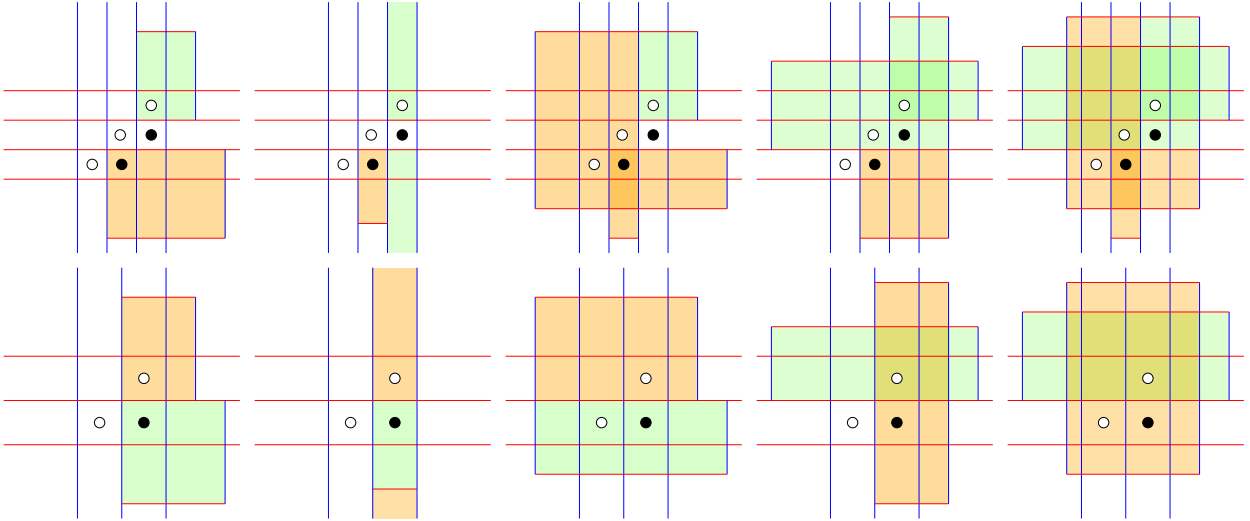}
\end{center}
\caption{The correspondence of domains in the first case. The domains
  appearing in the map $drs$ are shown in the top row and the
  corresponding domains appearing in $\Phi_{1,j}\Psi_{1,j}$ are shown
  in the bottom row. The domains corresponding to $s$ and $\Psi_{1,j}$
  are shown in green, while the domains corresponding to $d$ and
  $\Phi_{1,j}$ are shown in orange.}\label{fig:case1domains}
\end{figure}

Recall that the stabilization map $s$ comes from certain northeast
snail domains, and they are not allowed to pass through $z_{1,j}$;
similarly, the destabilization map $d$ comes from certain southeast
snail domains, and they too are not allowed to pass through
$z_{1,j}$. Therefore, only domains of certain `shapes' can appear in the
map $drs$, and these shapes correspond to the domains that appear in
the map $\Phi_{1,j}\Psi_{1,j}$. This is best illustrated by
\hyperref[fig:case1domains]{Figure \ref*{fig:case1domains}}.

\begin{figure}
\psfrag{wj}{$w_{1,j}$}
\psfrag{zj}{$z_{1,j}$}
\psfrag{wj1}{$w_{1,j+1}$}
\psfrag{zj1}{}
\psfrag{wja}{}
\psfrag{zja}{}
\psfrag{zjk}{$z_{1,k_1+1}$}
\psfrag{wjk}{$w_{1,k_1+1}$}
\psfrag{s}{$s$}
\psfrag{h}{$h$}
\psfrag{r}{$r$}
\psfrag{d}{$d$}
\psfrag{g1}{$\mc{G}^1$}
\psfrag{g2}{$\mc{G}^2$}
\psfrag{g3}{$\mc{G}^3=\mc{G}^4$}
\psfrag{g4}{$\mc{G}^5$}
\psfrag{g5}{$\mc{G}^6$}
\begin{center}
\includegraphics[width=0.9\textwidth]{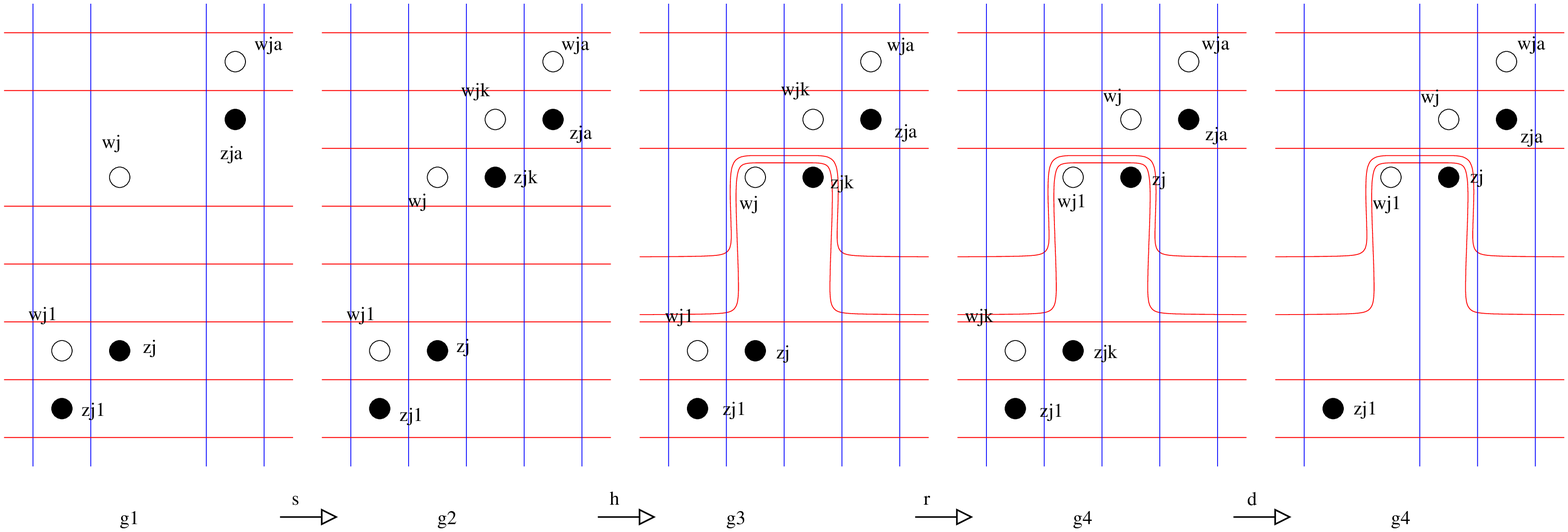}
\end{center}
\caption{The sequence of moves in the second case. There could be
arbitrarily many $\al$-circles between $w_{1,j}$ and $z_{1,j}$; we
have drawn three.}\label{fig:case2}
\end{figure}

\begin{figure}
\begin{center}
\includegraphics[width=0.8\textwidth]{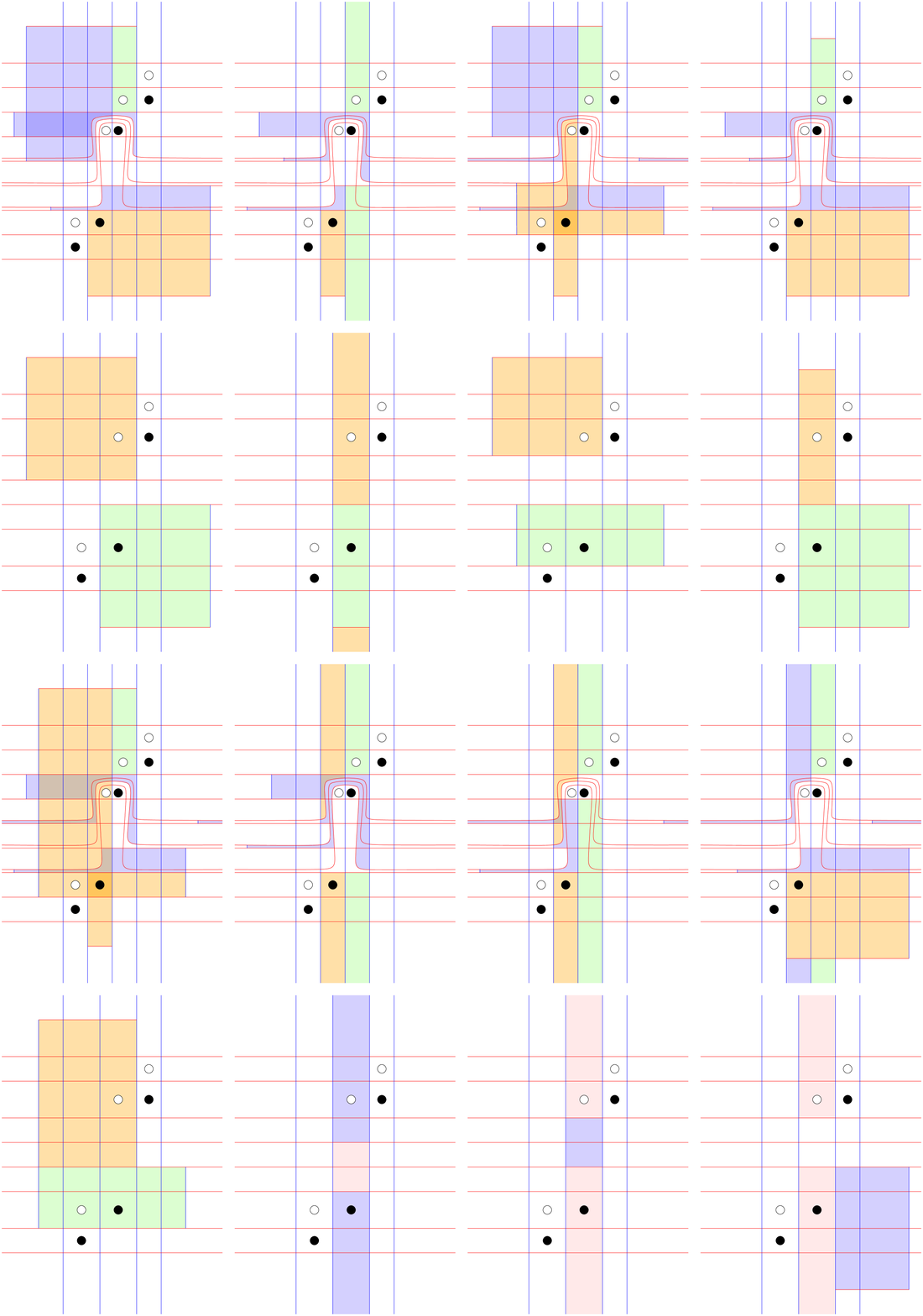}
\end{center}
\captionsetup{width=\textwidth}
\caption{The correspondence of domains in the second case (here we have
drawn four $\al$-circles between $w_{1,j}$ and $z_{1,j}$). The domains
  appearing in the map $drhs$ are shown in the first and the third row
  and the corresponding domains appearing in
  $\Phi_{1,j}\Psi_{1,j}+[H\colon\del]$ are shown in the second and the
  fourth row. The domains corresponding to $s$ and $\Psi_{1,j}$ are
  shown in green, the domains corresponding to $d$ and $\Phi_{1,j}$
  are shown in orange, the domains corresponding to $h$ and the first
  map in $[H\colon\del]$ are shown in violet, and the domains
  corresponding to the second map in $[H\colon\del]$ are shown in pink.}\label{fig:case2domains}
\end{figure}

In the second case, $\mc{G}^3=\mc{G}^4$ and $v=\Id$. The sequence of
moves is shown in \hyperref[fig:case2]{Figure \ref*{fig:case2}}. In
this case, the maps $\iota drhs$ and $c_{j,j+1}\Phi_{1,j}\Psi_{1,j}$
will not agree on the nose; we will have to modify the second map by a
homotopy. Given states $x,y\in\mc{S}_{\mc{G}_j}$, let $\mc{R}_H(x,y)$
be the subset of $\mc{R}_{\mc{G}_j}(x,y)$ consisting of all the
rectangles that contain both $w_{1,j}$ and $z_{1,j}$, with $z_{1,j}$
lying to the north of $w_{1,j}$. Define the
$U_{\imath,\jmath}$-equivariant map
$H\colon\CF_{\mc{G}_j}\rightarrow \CF_{\mc{G}_j}$ as follows:
$$H(x)=\sum_{y\in\mc{S}_{\mc{G}_j}}yU_{1,j}^{-1}\sum_{\substack{R\in\mc{R}_H(x,y)\\
n_z(R)=1}}\prod_{\imath,\jmath}U_{\imath,\jmath}^{n_{w_{\imath,\jmath}}(R)}.$$
The map $c_{j,j+1}H$ is an $U_i$-equivariant map from $\CF_{\mc{G}_j}$
to $\CF_{\mc{G}_{j+1}}$. Furthermore,
$[c_{j,j+1}H\colon\del]=(U_{1,j}+U_{1,j+1})c_{j,j+1} \Phi_{1,j}H +
c_{j,j+1}[H\colon\del]$. However, it is easy to see that
$\Phi_{1,j}H=0$ for the case at hand.  We will show that the maps
$\iota drhs$ and
$c_{j,j+1}\Phi_{1,j}\Psi_{1,j}+[c_{j,j+1}H\colon\del]=c_{j,j}(\Phi_{1,j}\Psi_{1,j}+[H\colon\del])$
are equal. The part $c_{j,j+1}$ once again comes from $r$ and
$d$. Therefore, we only have to show that the domains that contribute
to the map $drhs$ correspond to the domains that contribute to the map
$\Phi_{1,j}\Psi_{1,j}+[H\colon\del]$.

Recall that in a horizontal commutation, we change a single
$\al$-circle, and the commutation map is given by counting certain
pentagons, one of whose vertices is $\varrho$, which is one of the two
intersection points between the old $\al$-circle and the new
$\al$-circle.  Therefore, such pentagons could be two types:
the \emph{north pentagons} which lie to the north of $\varrho$, and
the \emph{south pentagons} which lie to the south of $\varrho$.  The
map $h$ is the composition of several such commutation
maps. Furthermore, it is clear that in the composition, a south
pentagon can not be followed by a north pentagon. Therefore, any
domain corresponding to the map $h$ is a sum of two domains: one
coming from the north pentagons, followed by one coming from the south
pentagons.

The north pentagons happen immediately after the stabilization map
$s$; hence there are only a few possible configurations for a domain
corresponding to the north pentagons. Similarly, the south pentagons
are immediately followed by the destabilization map $d$; therefore,
there are only a few possible configurations for a domain
corresponding to the south pentagons as well.  This allows us do a
case analysis, whereby we can show that the domains contributing to
map $drhs$ correspond to the domains contributing to the map
$\Phi_{1,j}\Psi_{1,j}+[H\colon\del]$,
cf.\ \hyperref[fig:case2domains]{Figure \ref*{fig:case2domains}}.
\end{proof}

To complete the proof of \hyperref[thm:main]{Theorem \ref*{thm:main}},
recall that the composition
$\wt{\tau}_{k_1,1}\cdots\wt{\tau}_{2,3}\wt{\tau}_{1,2}$ induces the
mapping class group element $\si_1\in\MCG(S^3,L,p)$. Therefore, the
composition $f_{k_1,1}\cdots f_{2,3}f_{1,2}$ must equal
$\rho(\si_1)$. However, by \hyperref[thm:changespecial]{Theorem
  \ref*{thm:changespecial}}, $f_{k_1,1}f_{k_1-1,k_1}\cdots f_{2,3}
f_{1,2}=f_{k_1,1}f_{k_1-1,k_1}\cdots
f_{1,3}=\cdots=f_{k_1,1}f_{1,k_1}=\Id+\Phi_1\Psi_1$.
\end{proof}

\section{Computations}\label{sec:computations}
Let us now present some computations. We will only work with knots in
$S^3$; since there is only one component to work with, we will write
$\Phi$ and $\Psi$ instead of $\Phi_1$ and $\Psi_1$. For a pointed knot
$(K,p)$, define the polynomials
$h_K(q,t),r_K(q,t)\in\Z[q,q^{-1},t,t^{-1}]$ as follows: the
coefficient of $q^at^b$ in $h_K(q,t)$ is the dimension (as an
$\F_2$-module) of $\wh{\HFL}(S^3,K,p)$ in $(M,A)$-bigrading $(a,b)$;
the coefficient of $q^at^b$ in $r_K(q,t)$ is the rank of the map
$\Phi\Psi$ on $\wh{\HFL}(S^3,K,p)$ in $(M,A)$-bigrading $(a,b)$. Since
$(\Phi\Psi)^2=0$, the only information that we can get from it is its
rank.

If a knot $K$ satisfies $h_K(q,t)=t^s\wt{h}(qt)$ for some $s\in\Z$ and
some $\wt{h}\in\Z[v,v^{-1}]$, it is called a \emph{thin knot}; if, in
addition, $s=-\frac{\si(K)}{2}$ where $\si(K)$ is the the signature of
$K$, then $K$ is called \emph{$\si$-thin}. The `Euler characteristic'
of knot Floer homology is its symmetrized Alexander-Conway polynomial
\cite{POZSzknotinvariants}, i.e.\ $h_K(-1,t)=\Delta_K(t)$. Therefore,
for $\si$-thin knots, the knot Floer homology can be reconstructed
from the Alexander-Conway polynomial and the knot signature as
$h_K(q,t)=(-q)^{\frac{\si(K)}{2}}\Delta_K(-qt)$. All quasi-alternating
knots are $\si$-thin \cite{POZSzalternatingknots,
  CMPOquasialternating}, and out of the $85$ prime knots up to nine
crossings, $83$ are quasi-alternating. The following theorem shows
that for thin knots, $r_K(q,t)$ can be constructed from $h_K(q,t)$,
and hence for $\si$-thin knots, $r_K(q,t)$ can be constructed from the
signature and the Alexander-Conway polynomial.

\begin{thm}\label{thm:computations}
  If $K$ is a thin knot with $h_K(q,t)=t^s\wt{h}(qt)$, then
  $r_K(q,t)=t^s\wt{r}(qt)$, where $\wt{r}\in\Z[v,v^{-1}]$ is the
  unique polynomial such that
  $\wt{h}(v)=\frac{1-v^{2m+1}}{1-v}+(v^{-1}+2+v)\wt{r}(v)$ for some
  $m\in\Z$.
\end{thm}

\begin{proof}
The uniqueness is the easy part. If there are two different
polynomials $\wt{r}$ and $\wt{r}'$ that satisfy the above equation,
then we get
$$(1+v)^2(1-v)(\wt{r}-\wt{r}')=v^{2m+2}-v^{2m'+2}.$$
Differentiating once and then putting $v=-1$ shows that $m=m'$ and
hence $\wt{r}=\wt{r}'$.

The knot Floer homology $\wh{\HFL}(S^3,K,p)$ is a graded object with
the grading being the Maslov grading\footnote{Since we are working
  with a thin knot, the Alexander grading differs from this grading by
  a constant.}, and due to \hyperref[lem:homotopyproperties]{Lemma
  \ref*{lem:homotopyproperties}}, it carries two commuting
differentials $\Phi$ and $\Psi$ of gradings $1$ and $-1$,
respectively. There are two knot-invariant spectral sequences: one
starts at $(\wh{\HFL}(S^3,K,p),\Phi)$ and converges to $\F_2$ where
the differentials on the $i\ith$ page shift the $(M,A)$-bigrading by
$(2i-1,i)$; and the second one starts at $(\wh{\HFL}(S^3,K,p),\Psi)$
and converges to $\F_2$ lying in Maslov grading $0$ where the
differentials on the $i\ith$ page shift the $(M,A)$-bigrading by
$(-1,-i)$.  Therefore, for thin knots, both the spectral sequences
must collapse immediately, i.e.\ $H_*(\wh{\HFL}(S^3,K,p),\Phi)=\F_2$
and $H_*(\wh{\HFL}(S^3,K,p),\Psi)=\F_2$.  Therefore, from a slight
generalization of \cite[Lemma 7]{IP}\footnote{\cite{IP} does it for
  bigraded complexes; however, the proof goes through for singly
  graded complexes like the one under consideration.}, we see that
$(\wh{\HFL}(S^3,K,p),\Phi,\Psi)$ must be isomorphic to a direct sum of
\emph{square} pieces of the form $(\F_2(a_1,a_2,a'_2,a_3),\Phi,\Psi)$,
where $\Phi(a_1)=\Psi(a_3)=a'_2$, $\Phi(a_2)=a_3$ and $\Psi(a_2)=a_1$,
and a single \emph{ladder} piece of the form
$(\F_2(b_1,\ldots,b_{2m+1}),\Phi,\Psi)$ where either
$\Phi(b_{2i-1})=\Psi(b_{2i+1})=b_{2i}$ for all $1\leq i\leq m$ or
$\Phi(b_{2i})=b_{2i+1},\Psi(b_{2i})=b_{2i-1}$ for all $1\leq i\leq
m$. The only contributions to the map $\Phi\Psi$ come from the square
pieces in the middle grading, viz.\ $\Phi\Psi(a_2)=a'_2$. This
produces the required decomposition of $\wt{h}(v)$: the square pieces
contribute the term $(v^{-1}+2+v)\wt{r}(v)$; the ladder piece
contributes the term $(1+\cdots+v^{2m})$ in the first case, and the
term $(1+\cdots+v^{-2m})$ in the second case.
\end{proof}

The two prime knots up to nine crossings that are not thin (and hence
not quasi-alternating) are $8_{19}$ and $9_{42}$.

\begin{lem}\label{lem:minorcompute}
  $r_{8_{19}}(q,t)=0$ and $r_{9_{42}}(q,t)=t^{-1}+q^2t.$
\end{lem}

\begin{proof}
From \cite{JBWG}, we get
\begin{align*}
h_{8_{19}}(q,t)&=t^{-3}+qt^{-2}+q^2+q^5t^2+q^6t^3,\\
h_{9_{42}}(q,t)&=q^{-1}t^{-2}+2t^{-1}+1+2q+2q^2t+q^3t^2.
\end{align*}

Since $(\Phi\Psi)^2=0$, the coefficient of $q^at^b$ in $r_K(q,t)$ is
less than or equal to half the coefficient of $q^at^b$ in
$h_K(q,t)$. This immediately shows that $r_{8_{19}}(q,t)=0$.

The knot Floer homology of $9_{42}$ is supported on two diagonals,
$M-A=0$ and $M-A=1$. For $i\in\{0,1\}$, let $C_i$ be the direct
summand that is supported in the diagonal $M-A=i$. As in the proof of
\hyperref[thm:computations]{Theorem \ref*{thm:computations}}, there
are two knot-invariant spectral sequences starting at
$(C_1,\Phi)\oplus(C_0,\Phi)$ and $(C_1,\Psi)\oplus(C_0,\Psi)$ and
ending at $\F_2$. \cite{JBWG} tells us that the first spectral
sequence collapses at the very next step, i.e.\
$H_*(C_1,\Phi)=0$. Since the knot $9_{42}$ is reversible and the
spectral sequences are knot invariants, the second spectral sequence
must also collapse at the very next step, i.e.\
$H_*(C_1,\Psi)=0$. Another application of \cite[Lemma 7]{IP} tells us
that $(C_1,\Phi,\Psi)$ is a direct sum of two square pieces, and hence
$r_{9_{42}}(q,t)=t^{-1}+q^2t.$
\end{proof}

We conclude with the following observation.
\begin{thm}
  There does not exist any orientation preserving involution of $S^3$
  that acts freely and in an orientation preserving way on either
  the knot $8_{20}$ or the knot $9_{42}$.
\end{thm}

\begin{proof}
  Let $(K,p)$ be a pointed oriented knot in $S^3$, and let $\tau$ be
  an orientation preserving involution of $S^3$ that acts freely and
  in an orientation preserving way on $K$. Let $\ol{\si}$ be half a
  Dehn twist around $K$ that interchanges the points $p$ and $\tau(p)$
  such that $\ol{\si}^2=\si$, the full positive Dehn twist around
  $K$. Since $\tau$ acts freely on a neighborhood of $K$ and since
  $\ol{\si}$ is identity outside a neighborhood of $K$, the two
  diffeomorphisms $\tau\ol{\si}$ and $\ol{\si}\tau$ induce the same
  element of $\MCG(S^3,K,\tau(p))$; hence, the two diffeomorphisms
  $(\tau\ol{\si})^2$ and $\tau^2\ol{\si}^2=\si$ induce the same
  element in $\MCG(S^3,K,p)$.

  Let
  $\wh{\rho}\colon\MCG(S^3,K,p)\rightarrow\Aut_{N(\mc{C}_1)}(\wh{\HFL}(S^3,K,p))$
  be the action of the mapping class group on knot Floer
  homology. \hyperref[thm:main]{Theorem \ref*{thm:main}} tells us that
  $\wh{\rho}(\si)=\Id+\Phi\Psi$. However, since $(\tau\ol{\si})^2=\si$
  in the mapping class group,
  $\wh{\rho}(\si)=\Id+\Phi\Psi=(\wh{\rho}(\tau\ol{\si}))^2$.

  We know from \cite{JBWG}, \hyperref[thm:computations]{Theorem
    \ref*{thm:computations}} and \hyperref[lem:minorcompute]{Lemma
    \ref*{lem:minorcompute}},
\begin{align*}
h_{8_{20}}(q,t)&=q^{-2}t^{-2}+2q^{-1}t^{-1}+3+2qt+q^2t^2\\
r_{8_{20}}(q,t)&=q^{-1}t^{-1}+qt\\
h_{9_{42}}(q,t)&=q^{-1}t^{-2}+2t^{-1}+1+2q+2q^2t+q^3t^2\\
r_{9_{42}}(q,t)&=t^{-1}+q^2t.
\end{align*}

For the knot $8_{20}$ in $(M,A)$-bigrading $(1,1)$, and for the knot
$9_{42}$ in $(M,A)$-bigrading $(2,1)$, the homology is two dimensional
and the rank of the map $\Phi\Psi$ is one. Therefore, in either case,
we can do a change of basis to represent the map $\Id+\Phi\Psi$ by the
matrix
$$
\begin{pmatrix}
1&1\\
0&1
\end{pmatrix}
.$$
However, this matrix is not a square in $\GL_2(\F_2)$.
\end{proof}

\bibliographystyle{amsalpha}
\bibliography{action}

\end{document}